\documentclass[]{article}
\usepackage{ctex}
\usepackage{amsmath,amssymb}
\usepackage{geometry}
\usepackage{xcolor}
\definecolor{deepgreen}{RGB}{0,100,0}
\usepackage{changepage}
\usepackage{quoting}
\usepackage{enumitem}
\usepackage{setspace}
\usepackage{relsize}
\usepackage{extarrows}

\numberwithin{equation}{section}
\usepackage{tikz}
\usepackage{ctex}
\usepackage{graphicx}
\usepackage{pifont}
\usepackage{amssymb}
\usepackage{booktabs}
\usepackage{algorithm}
\usepackage{algorithmic}
\usepackage{graphicx}
\graphicspath{{./paper-picture/}}
\usepackage{multirow}  
\usepackage{tikz}
\definecolor{mgq}{rgb}{0,0,0}
\newcommand{\mgq}{\color{mgq}}
\usetikzlibrary{decorations.pathmorphing}

\title{{\mgq Fourier--U-Net: a hybrid method for acoustic scattering media reconstruction}}
\author{Wangling Zhang\footnotemark[2], Guanqiu Ma\footnotemark[1]\; \footnotemark[2]}
	
	\footnotetext[1] {Corresponding author}
	\footnotetext[2]{School of Mathematical Sciences, Sichuan Normal University, Chengdu, 610068, China(wl\_zhang@stu.sicnu.edu.cn, gqma@sicnu.edu.cn).}
	
\date{}

\begin{document}
	\maketitle
	\begin{abstract}
		We {\mgq propose} a Fourier--U-Net {\mgq hybrid} method for {\mgq acoustic scattering media reconstruction}. 
		Under the Born approximation, {\mgq an initial contrast estimate is constructed via a truncated Fourier series: the nonzero-frequency Fourier coefficients are estimated from far-field scattering data, while the zero-frequency coefficient is supplied by a low-frequency compensation term derived from the far-field value at a small wave number, with finite-distance measurement effects on far-field extraction accounted for.} A U-Net then learns the mapping from this initial estimate to the true contrast {\mgq suppressing blurring, ringing artifacts, and contrast bias caused by Fourier truncation, finite-distance far-field approximation, and measurement noise. Numerical examples validate the method in terms of parameter selection under different receiver distances and truncation orders, adaptability under noise mismatch, reconstruction of complex targets, and generalization to unseen categories, indicating its effectiveness and feasibility.}
		
		%Four numerical experiments evaluate the method. The effects of receiver distance and Fourier truncation order on far-field extraction errors, reconstruction accuracy, and the number of Fourier modes are examined, together with the role of zero-frequency compensation through an ablation test. High-noise pretraining followed by fine-tuning at the target noise level is also evaluated. MNIST handwritten digits are used to test the reconstruction of complex targets. The MNIST-trained model is then tested on EMNIST handwritten letters and Chinese handwritten numerals without further training to assess transfer to unseen character categories. For synthetic scattering data generated under the Born approximation, with finite-distance measurements and fixed relative noise levels, the method reduces reconstruction errors compared with direct Fourier reconstruction and improves structural similarity for the character targets.

	\end{abstract}

	\section{Introduction}
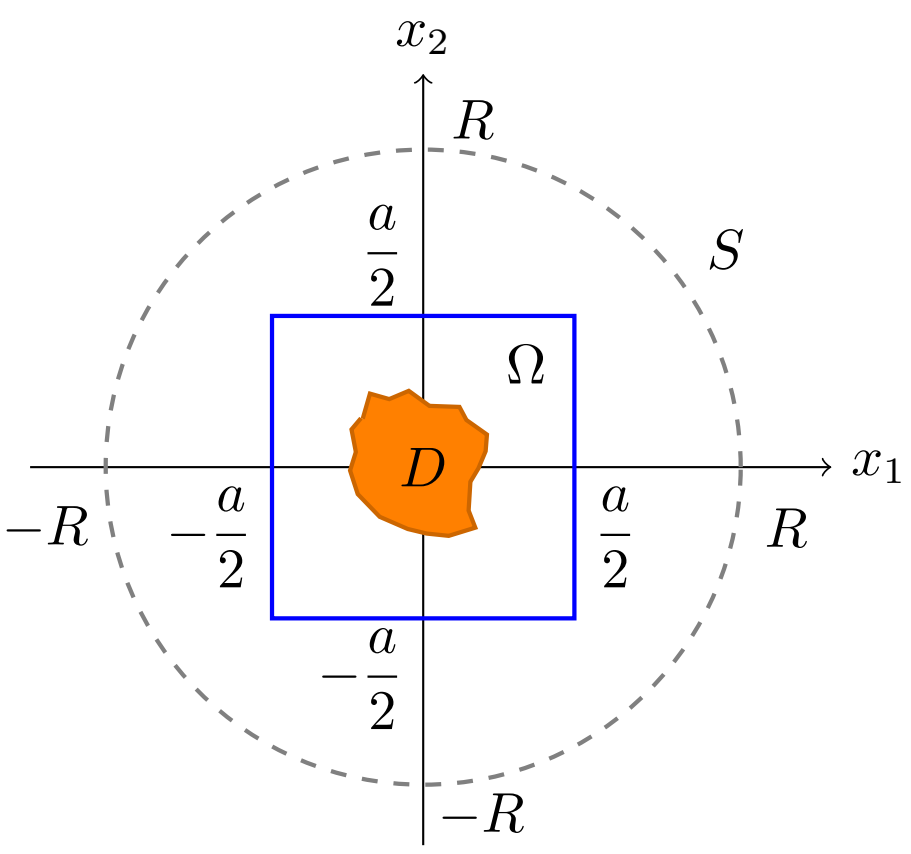
\begin{figure}[htbp]
		\centering
		
		\begin{tikzpicture}[scale=2]
			\draw[->] (-1.3,0) -- (1.35,0) node[right] {$x_1$};
			\draw[->] (0,-1.25) -- (0,1.3) node[above] {$x_2$};
			
			\draw[dashed,gray,thick] (0,0) circle (1.05);
			
			\draw[blue,line width=0.8pt]
			(-0.5,-0.5) rectangle (0.5,0.5);
			
			\node[left] at (0,0.7) {$\dfrac{a}{2}$};
			\node[left] at (0,-0.7) {$-\dfrac{a}{2}$};
			
			\node[below left] at (-0.5,0) {$-\dfrac{a}{2}$};
			\node[below right] at (0.5,0) {$\dfrac{a}{2}$};
			
			\node at (0.34,0.34) {$\Omega$};
			
			\filldraw
			[
			fill=orange,
			draw=orange!80!black,
			thick,
			decorate,
			decoration={
				random steps,
				segment length=4pt,
				amplitude=1.5pt
			}
			]
			plot[smooth cycle] coordinates{
				(-0.22,0.05)
				(-0.20,0.16)
				(-0.11,0.22)
				(0.02,0.20)
				(0.14,0.15)
				(0.20,0.05)
				(0.16,-0.05)
				(0.18,-0.14)
				(0.08,-0.22)
				(-0.05,-0.20)
				(-0.13,-0.14)
				(-0.22,-0.10)
			};
			
			\node at (0,0) {$D$};
			\node at (1.0,0.72) {$S$};
			\node at (1.20,-0.2) {$R$};
			\node at (-1.25,-0.2) {$-R$};
			\node[left] at (0.3,1.15) {$R$};
			\node[left] at (0.4,-1.15) {$-R$};
			
		\end{tikzpicture}
		
		\caption{Schematic illustration of the two-dimensional imaging domain $\Omega$, scatterer $D$, and measurement circle $S$.}
		\label{fig:measurement_geometry}
	\end{figure}

	Inverse medium scattering aims to recover unknown medium parameters from known incident fields and measured scattered fields. It has applications in electromagnetics, optics, medical imaging, and geophysical exploration.
	As shown in Fig.~\ref{fig:measurement_geometry}, {\mgq we consider acoustic scattering in a homogeneous background medium occupying $\mathbb{R}^2$. Let $\Omega\subset\mathbb{R}^2$ be a bounded square imaging domain and let $D\subset\Omega$ denote an inhomogeneous scatterer whose material properties differ from those of the background. The scatterer is characterized by the refractive index function $\varepsilon_r(\mathbf{x})$, which satisfies
	\[
	\varepsilon_r(\mathbf{x})\geq1,\qquad
	\mathrm{supp}\!\left(\varepsilon_r(\mathbf{x})-1\right)\subset\Omega,
	\]
	so that the scatterer $D$ lies entirely within $\Omega$. Outside $D$, the medium is homogeneous with $\varepsilon_r=1$.}
	The contrast function is defined by
	$\chi(\mathbf{x})=\varepsilon_r(\mathbf{x})-1$,
	which represents the {\mgq refractive index} difference between the scatterer and the background.
	The incident field is a time-harmonic plane wave
	$u^{i}(\mathbf{x})=e^{ik\hat{\theta}\cdot\mathbf{x}}$,
	where $k$ is the {\mgq wave number} and
	{\mgq$\hat{\theta}\in\mathbb{S}$} is the incident direction.
	
	The total field $u(\mathbf{x})$ is the sum of the incident and scattered fields,
	$u(\mathbf{x})=u^{i}(\mathbf{x})+u^{s}(\mathbf{x})$,
	and satisfies the Helmholtz equation
	\begin{equation}
		\Delta u(\mathbf{x})+k^2\varepsilon_r(\mathbf{x})u(\mathbf{x})=0,
		\qquad
		\mathbf{x}\in\mathbb{R}^2.
		\label{eq:helmholtz}
	\end{equation}
	The scattered field $u^{s}(\mathbf{x})$ satisfies the Sommerfeld radiation condition
	\begin{equation}
		\lim_{r=|\mathbf{x}|\rightarrow\infty}
		\sqrt{r}
		\left(
		\frac{\partial u^s}{\partial r}-iku^s
		\right)=0.
		\label{eq:sommerfeld}
	\end{equation}
	The scattering problem can equivalently be written as the Lippmann--Schwinger integral equation
	\begin{equation}
		u(\mathbf{x})
		=
		u^{i}(\mathbf{x})
		+
		k^2\int_{\Omega}
		G(\mathbf{x},\mathbf{y})
		\chi(\mathbf{y})
		u(\mathbf{y})\,d\mathbf{y},
		\qquad
		\mathbf{x}\in\mathbb{R}^2,
		\label{eq:ls}
	\end{equation}
	where $G(\mathbf{x},\mathbf{y})$ is {\mgq the Green's function}
	\begin{equation}
		G(\mathbf{x},\mathbf{y})
		=
		\frac{i}{4}H_0^{(1)}(k|\mathbf{x}-\mathbf{y}|),
		\label{eq:green}
	\end{equation}
	and $H_0^{(1)}$ denotes the zeroth-order Hankel function of the first kind.
	On the measurement circle $S$, the scattered field is expressed as
	\begin{equation}
		u^{s}(\mathbf{x})
		=
		k^2\int_{\Omega}
		G(\mathbf{x},\mathbf{y})
		\chi(\mathbf{y})
		u(\mathbf{y})\,d\mathbf{y},
		\qquad
		\mathbf{x}\in S.
		\label{eq:scattered_field}
	\end{equation}
	
	Let $a$ denote the side length of the imaging domain, so that
	$\Omega=[-\frac{a}{2},\frac{a}{2}]\times[-\frac{a}{2},\frac{a}{2}]$.
	The receivers are located on a circle centered at the origin with radius $R$.
	For an observation direction {\mgq $\hat{d}\in\mathbb{S}$}, the receiver position is
	\begin{equation}
		\mathbf{x}=R\hat{d},\qquad R=|\mathbf{x}|.
	\end{equation}
	The measurement circle enclosing $\Omega$ is defined by
	\begin{equation}
		S=
		\left\{
		\mathbf{x}=R\hat{d}:
		\hat{d}\in\mathbb{S}
		\right\},
		\qquad
		\overline{\Omega}\subset B_R(0).
		\label{eq:measurement_surface}
	\end{equation}
	Suppose that $N_{\mathrm{inc}}$ incident-wave configurations are used.
	For the $i$-th configuration, let $k_i$ and $\hat{\theta}_i$ denote the wave number and incident direction, respectively.
	The corresponding scattered-field data on $S$ are given by
	\begin{equation}
		\mathcal{D}
		=
		\left\{
		u^s\!\left(\mathbf{x};k_i,\hat{\theta}_i\right)
		\,\middle|\,
		\mathbf{x}\in S,\;
		i=1,2,\ldots,N_{\mathrm{inc}}
		\right\}.
		\label{eq:data}
	\end{equation}
	The inverse medium problem (\textbf{IMP}) considered in this paper is as follows:
	
	Given the incident fields
	$\left\{
	u^{i}(\mathbf{x};k_j,\hat{\theta}_j)
	\right\}_{j=1}^{N_{\mathrm{inc}}}$
	{\mgq with the finite} wave numbers and incident directions,
	and the corresponding scattered-field data
	$\left\{
	u^{s}(\mathbf{x};k_j,\hat{\theta}_j)
	\right\}_{j=1}^{N_{\mathrm{inc}}}$
	on the measurement circle $S$,
	the objective is to estimate the medium contrast field
	$\chi(\mathbf{x})$ within $\Omega$ and identify its support.

	{\mgq The scattered field depends non-linearly on the medium contrast; moreover, the inverse medium scattering problem is ill-posed}~\cite{colton2019inverse}.
	Measurement errors and noise, together with limited incident directions, observation directions, and frequencies, restrict the information available for reconstruction.
	Multi-frequency measurements provide additional scattering information and can improve {\mgq the stability and quality of the reconstruction}~\cite{bao2015multifrequency}.
	Accurate recovery of medium parameters from finite and noisy scattering data nevertheless remains challenging.
	
	Inverse scattering reconstruction methods are commonly divided into iterative and non-iterative methods.
	Iterative methods repeatedly solve the forward problem and update the medium parameters to reduce the discrepancy between computed and measured scattered fields.
	These methods can account for nonlinear scattering, but often require substantial computation and depend on the initialization and regularization strategy~\cite{colton2019inverse,vandenberg1997contrast}.
	Non-iterative methods relate the measured data directly to target information through analytical relations or suitable model approximations, avoiding repeated parameter updates.
	Direct sampling methods, for example, use indicator functions to locate scatterers efficiently and can tolerate measurement noise~\cite{potthast2006survey,ito2012sampling,ito2013electromagnetic,li2013farfield}.
	Fourier-based methods use the relation between scattering data and target spatial frequencies to estimate Fourier coefficients and construct an initial reconstruction.
	The resulting procedure {\mgq entails a low computational cost} and has a direct interpretation in terms of spatial frequencies.
	Fourier methods have been developed for reconstructing unknown sources from multi-frequency data~\cite{zhang2015fourier,wang2017fourier}.
	The windowed Fourier transform has also been used for source reconstruction from single-frequency far-field data~\cite{griesmaier2012windowed}.
	Multi-frequency data can improve the stability of inverse source problems~\cite{cheng2016stability}.
	For diffraction tomography, Devaney developed a filtered backpropagation algorithm under the first Born or Rytov approximation~\cite{devaney1982filtered}.
	Limited incident and observation angles and a finite frequency range restrict the spatial-frequency information that can be recovered.
	The numerical examples below illustrate that finite-mode reconstructions may {\mgq suffer} contrast errors, blurred boundaries, and ringing artifacts.
	With limited scattering data, compensating for missing frequency information and suppressing artifacts remain challenges for accurate contrast reconstruction.
	
	Neural networks can learn nonlinear mappings from approximate reconstructions to true contrasts and use the prior information contained in training data to improve reconstruction quality~\cite{chen2020review,jin2017cnn,mccann2017review,arridge2019data}.
	In inverse scattering, deep learning methods have been applied to reconstruction with nonlinear scattering models, limited-view data, and noisy measurements~\cite{wei2019deep,yao2019twostep,zhou2021improved,gao2022ann}.
	U-Net uses an encoder--decoder architecture with skip connections to combine features at different spatial scales, allowing both global structures and local details to be represented~\cite{ronneberger2015unet}.
	This architecture is adopted here to refine Fourier initial reconstructions affected by blurring, ringing, and contrast errors.
	
	{\mgq In this paper, the Fourier--U-Net method is proposed for acoustic scattering media reconstruction. Under the Born approximation, we construct an initial contrast estimate through a truncated Fourier series. The nonzero Fourier coefficients are estimated from far-field patterns, with the incident and observation configurations chosen according to the retained nonzero Fourier modes, whereas the zero-frequency coefficient is approximated from the far-field value at a small wave number.}
	The real and imaginary parts of the normalized initial estimate form the two input channels of a U-Net.
	The network predicts the normalized target contrast, and its output is rescaled to recover the contrast.
	Benny et al. used the Fourier diffraction theorem to construct an initial contrast estimate and a U-Net to refine microwave inverse scattering reconstructions~\cite{benny2022microwave}.
	Other studies have combined initial images from physical models with neural networks for quantitative inverse scattering~\cite{zhang2020hybrid,ning2024dsm}.
	For the multi-frequency inverse source problem, Chen et al. combined Fourier initial reconstruction and U-Net refinement with high-noise pretraining followed by fine-tuning at the target noise level~\cite{chen2026fourier}.
	{\mgq Throughout this paper, all scattered-field data used in the numerical experiments are synthesized under the Born approximation. The method is therefore assessed on Born-synthetic data with finite receiver distance, Fourier truncation, and measurement noise. Its applicability to strongly scattering media, where the Born approximation may not hold, is discussed as a limitation in Section 5.}
	
	The remainder of this paper is organized as follows.
	The Fourier initial reconstruction method is introduced in Section~2.
	The U-Net architecture and training strategy for noisy data are presented in Section~3.
	Numerical experiments are reported and their results are analyzed in Section~4.
	The conclusions, limitations, and future work are discussed in Section~5.

	\section{Fourier Initial Reconstruction}
	The inverse scattering model is linearized using the Born approximation to obtain an initial estimate of the target contrast.
	The relation between far-field scattering patterns and the Fourier transform is used to estimate a finite set of Fourier coefficients.
	A truncated Fourier series constructed from these coefficients provides the initial reconstruction for subsequent U-Net refinement.
	
	Under the Born approximation, the total field in the scattering region is approximated by the incident field,
	$u(\mathbf{y})\approx u^{i}(\mathbf{y})$,
	so that the scattered field is approximated by
	\begin{equation}
		u^s(\mathbf{x})\approx k^2 \int_{\Omega} G(\mathbf{x},\mathbf{y})\chi(\mathbf{y})
		e^{ik\hat{\theta}\cdot\mathbf{y}}\,d\mathbf{y},\qquad\mathbf{x}\in S.
		\label{eq:born}
	\end{equation}
	The approximate scattered field depends linearly on the target contrast, reducing the reconstruction of $\chi(\mathbf{x})$ from the measured scattered fields to a linear inverse problem.
	When $k|\mathbf{x}-\mathbf{y}|\rightarrow\infty$, the two-dimensional Green's function has the asymptotic form
	\begin{equation}
		G(\mathbf{x},\mathbf{y})\sim\frac{e^{i\pi/4}}{\sqrt{8\pi k}}
		\frac{e^{ik|\mathbf{x}-\mathbf{y}|}}
		{\sqrt{|\mathbf{x}-\mathbf{y}|}}.
		\label{eq:green_asym}
	\end{equation}
	Substituting Eq.~\eqref{eq:green_asym} into Eq.~\eqref{eq:born} and taking the far-field limit $R\rightarrow\infty$ gives
	\begin{equation}
		u^s(\mathbf{x})\approx\frac{1}{\sqrt{R}}e^{ikR}E_\infty(\hat{d};\hat{\theta},k),
		\label{eq:scattered_farfield}
	\end{equation}
	where $E_\infty(\hat{d};\hat{\theta},k)$ denotes the far-field pattern under the Born approximation,
	\begin{equation}
		E_\infty(\hat{d};\hat{\theta},k)=C(k)\int_{\Omega}\chi(\mathbf{y})e^{-ik(\hat{d}-\hat{\theta})\cdot\mathbf{y}}\,d\mathbf{y},
		\label{eq:farfield_pattern}
	\end{equation}
	with $C(k)=\frac{k^{3/2}e^{i\pi/4}}{\sqrt{8\pi}}$.
	
	Equation~\eqref{eq:farfield_pattern} shows that the far-field pattern is proportional to the Fourier transform of the target contrast.
	Discrete sampling over a finite spatial-frequency range is therefore achieved by selecting the incident direction $\hat{\theta}$, observation direction $\hat{d}$, and wave number $k$.
	
	Let $\mathbf{x}=(x_1,x_2)$ denote the spatial coordinates.
	For each nonzero Fourier mode index
	$l=(l_x,l_y)\in\mathbb{Z}^2\setminus\{(0,0)\}$,
	the measurement configuration is defined by
	\begin{equation}
		\hat{\theta}_l = -\frac{l}{|l|},\qquad
		\hat{d}_l = \frac{l}{|l|},\qquad
		k_l=\frac{\pi|l|}{a}.
		\label{eq:measurement_config}
	\end{equation}
	This configuration satisfies
	\[
	k_l(\hat{d}_l-\hat{\theta}_l)=\frac{2\pi l}{a}.
	\]
	Equation~\eqref{eq:farfield_pattern} then becomes
	\begin{equation}
		E_\infty(\hat{d}_l;\hat{\theta}_l;k_l)=C(k_l)\int_{\Omega}\chi(\mathbf{y})
		e^{-i\frac{2\pi}{a}l\cdot\mathbf{y}}\,d\mathbf{y}.
		\label{eq:E_inf_fourier}
	\end{equation}
	The corresponding Fourier coefficient of the target contrast is defined by
	\begin{equation}
		\hat{\chi}_l=\frac{1}{a^2}\int_{\Omega}\chi(\mathbf{y})e^{-i\frac{2\pi}{a}l\cdot\mathbf{y}}\,d\mathbf{y},
		\label{eq:fourier_definition}
	\end{equation}
	which gives, for all nonzero modes $l\neq(0,0)$,
	\begin{equation}
		\hat{\chi}_l=\frac{E_\infty(\hat{d}_l;\hat{\theta}_l;k_l)}{a^2C(k_l)}.
		\label{eq:fourier_coeff}
	\end{equation}
	
	For the zero-frequency mode $l=(0,0)$, the direction $l/|l|$ is undefined, so Eq.~\eqref{eq:measurement_config} cannot be applied.
	The zero-frequency coefficient is approximated {\mgq by the} far-field value at a small wave number.
	A sufficiently small nonzero wave number $k_0$ and a reference direction $\hat{\theta}_0$ are selected, with the observation direction set to $\hat{d}_0=-\hat{\theta}_0$.
	This gives the approximation
	\begin{equation}
		\widehat{\chi}_{\boldsymbol{0}}\approx\frac{E_\infty(-\hat{\theta}_0;\hat{\theta}_0;k_0)}{a^2C(k_0)}.
		\label{eq:zero_frequency_approximation}
	\end{equation}
	The right-hand side is the Fourier transform of $\chi$ evaluated at the spatial frequency $-2k_0\hat{\theta}_0$ and divided by $a^2$.
	For sufficiently small $k_0$, this value approximates the zero-frequency Fourier coefficient.
	{\mgq The far-field value }at $k_0$ is denoted by
	\begin{equation}
		E_0=E_\infty(-\hat{\theta}_0;\hat{\theta}_0;k_0).
		\label{eq:zero_frequency_farfield}
	\end{equation}
	The zero-frequency contribution is treated as an independent low-frequency compensation term derived from $E_0$.
	It is not obtained through finite-distance far-field extraction.
	The subsequent finite-distance error assessment therefore includes only nonzero Fourier modes.
	
	For finite-distance measurements, the scattered-field asymptotic relation gives the approximation
	\begin{equation}
		E_\infty(\hat{d};\hat{\theta},k)\approx u^s(\mathbf{x})\sqrt{R}e^{-ikR}.
		\label{eq:extract_farfield}
	\end{equation}
	The initial contrast estimate is then formed from the estimated coefficients within
	$|l_x|\leq N$ and $|l_y|\leq N$
	using the truncated Fourier series
	\begin{equation}
		\chi_{\mathrm{rec}}(\mathbf{x})=\sum_{l_x=-N}^{N}\sum_{l_y=-N}^{N}\hat{\chi}_{l_x,l_y}e^{i\frac{2\pi}{a}(l_xx_1+l_yx_2)}.
		\label{eq:ifft}
	\end{equation}
	The truncation order $N$ gives a total of $(2N+1)^2$ Fourier coefficients.
	The zero-frequency coefficient is supplied by the theoretical low-frequency compensation term, while the remaining $(2N+1)^2-1$ coefficients are estimated from finite-distance scattering data.
	Increasing $N$ retains more high-frequency information and can improve spatial resolution, but requires more nonzero-mode measurements and increases the computational cost.
	Retaining additional high-frequency terms may also increase the influence of noise and finite-distance far-field approximation errors.
	The truncation order should therefore balance reconstruction accuracy, {\mgq computational cost}, and data stability.
	These effects are examined in the numerical experiments.

	\section{U-Net Network and Noisy-Data Training Strategy}
	The preceding section constructs an initial contrast estimate from a finite set of Fourier coefficients under the Born approximation.
	{\mgq The} Fourier truncation, finite-distance far-field approximation, and measurement noise can produce blurred boundaries, ringing artifacts, and contrast bias in this estimate.
	A U-Net is therefore used to learn the mapping from the Fourier initial reconstruction to the true target contrast.
	
	This section describes the U-Net input representation, architecture, loss function, and training procedure.
	Training strategies for matching and differing training and testing noise levels are then discussed.
	
	\subsection{U-Net Architecture and Training Procedure}
	
	The Fourier initial reconstruction may have a nonzero imaginary part, so both its real and imaginary parts are retained as separate input channels.
	The initial reconstruction and the true target contrast are normalized by the same scale
	$\chi_{\mathrm{s}}$:
	\begin{equation}
		\widetilde{\chi}_{\mathrm{rec}}=\frac{\chi_{\mathrm{rec}}}{\chi_{\mathrm{s}}},\qquad \widetilde{\chi}_{\mathrm{true}}=\frac{\chi_{\mathrm{true}}}{\chi_{\mathrm{s}}}.
		\label{eq:normalization}
	\end{equation}
	Here, $\chi_{\mathrm{s}}>0$ is the normalization scale, whose value is specified in the numerical experiments.
	The network input is written as
	\begin{equation}
		\mathbf{x}=\left[\operatorname{Re}\{\widetilde{\chi}_{\mathrm{rec}}\},\operatorname{Im}\{\widetilde{\chi}_{\mathrm{rec}}\}\right],
		\label{eq:unet_input}
	\end{equation}
	where the two components form the input channels of the network.
	The nonlinear mapping learned by the U-Net is represented by
	\begin{equation}
		\widetilde{\chi}_{\mathrm{final}}=\mathcal{U}_{\Theta}(\mathbf{x}),
		\label{eq:unet_mapping}
	\end{equation}
	where $\mathcal{U}_{\Theta}$ denotes the U-Net with parameters $\Theta$.
	
	\begin{figure}[htbp]
		\centering
		\includegraphics[width=\textwidth]{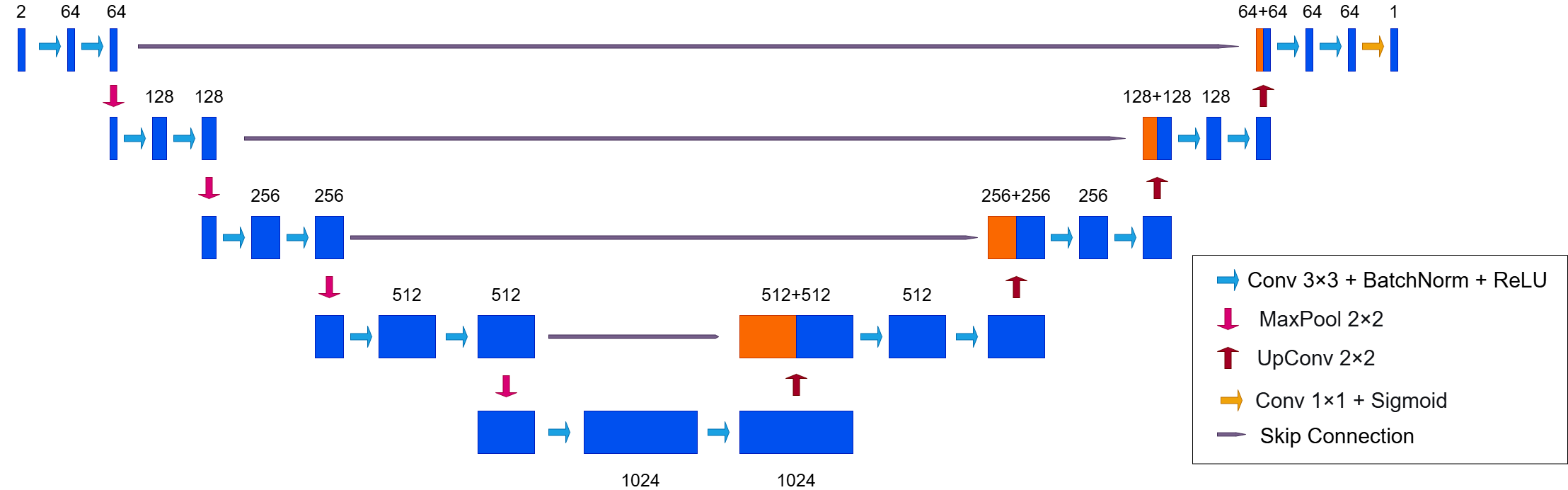}
		\caption{U-Net architecture}
		\label{fig:unet_visual}
	\end{figure}
	
	As shown in Fig.~\ref{fig:unet_visual}, the network uses an encoder--decoder architecture~\cite{ronneberger2015unet}.
	It has $2$ input channels, with $64$, $128$, $256$, and $512$ feature channels in the encoder and $1024$ channels at the bottleneck.
	Each convolutional block contains two $3\times3$ convolutions, each followed by batch normalization (BN) and a ReLU activation.
	The $3\times3$ convolutions use no bias terms.
	The encoder downsamples the feature maps using $2\times2$ max pooling.
	The decoder upsamples them using $2\times2$ transposed convolutions and concatenates them with the corresponding encoder feature maps through skip connections.
	The decoder channel count decreases to $64$, and a $1\times1$ convolution produces a single-channel output.
	
	A Sigmoid activation is applied at the network output to match the range of the normalized target labels:
	\begin{equation}
		\widetilde{\chi}_{\mathrm{final}}=\sigma(f_{\Theta}(\mathbf{x})),\qquad\sigma(z)=\frac{1}{1+e^{-z}},
		\label{eq:sigmoid_output}
	\end{equation}
	where $f_{\Theta}$ denotes the network output before Sigmoid activation.
	The target contrast is recovered at its original scale by
	\begin{equation}
		\chi_{\mathrm{final}}=\chi_{\mathrm{s}}\widetilde{\chi}_{\mathrm{final}}.
		\label{eq:output_rescale}
	\end{equation}

	The U-Net is trained by minimizing the mean squared error (MSE) between the predicted and true normalized contrasts.
	For a reconstruction region containing $N_p$ pixels, the loss is defined as
	
	\begin{equation}
		\mathcal{L}_{\mathrm{MSE}}=\frac{1}{N_p}\sum_{j=1}^{N_p}\left|\widetilde{\chi}_{\mathrm{final},j}-\widetilde{\chi}_{\mathrm{true},j}\right|^2.
		\label{eq:mse_loss}
	\end{equation}
	The total training loss is
	$\mathcal{L}_{\mathrm{total}}=\mathcal{L}_{\mathrm{MSE}}$.
	The corresponding optimization problem is written as
	\begin{equation}
		\Theta^{*}=\arg\min_{\Theta}\mathcal{L}_{\mathrm{total}}(\Theta).
		\label{eq:optimization}
	\end{equation}
	
	The network parameters are updated using the Adam optimizer.
	At update $t$, the parameter update is written as
	\begin{equation}
		\Theta^{(t+1)}
		\leftarrow
		\Theta^{(t)}
		-
		\eta_t
		\frac{\widehat{m}_t}
		{\sqrt{\widehat{v}_t}+\epsilon_{\mathrm{Adam}}}.
		\label{eq:parameter_update}
	\end{equation}
	Here, $\widehat{m}_t$ and $\widehat{v}_t$ are the bias-corrected first and second moment estimates of the accumulated gradients.
	The learning rate is denoted by $\eta_t$, and $\epsilon_{\mathrm{Adam}}>0$ is a numerical stabilizer.
	The square root and division in this expression are applied componentwise.
	
	The learning rate is adjusted using a step decay schedule.
	At epoch $e$, it is given by
	\begin{equation}
		\eta_e=\eta_0\gamma^{\left\lfloor (e-1)/T \right\rfloor},\qquad e=1,2,\ldots,
		\label{eq:lr_schedule}
	\end{equation}
	where $\eta_0$ is the initial learning rate,
	$T$ is the decay period in epochs,
	$\gamma\in(0,1)$ is the decay factor,
	and $\lfloor\cdot\rfloor$ denotes the floor function.
	The training parameters are specified in the numerical experiments.
	
	An independent validation set is used to monitor network performance during training.
	The parameters with the lowest validation loss are retained as the final model.
	When enabled, early stopping terminates training if the validation loss does not decrease for a prescribed number of consecutive epochs.
	
	The training procedure is summarized in Algorithm~\ref{alg:training_flow}.
	\begin{algorithm}[htbp]
		\caption{U-Net training procedure}
		\label{alg:training_flow}
		\begin{algorithmic}[1]
			\STATE \textbf{Input:} $\mathcal{D}_{\mathrm{train}}$ and $\mathcal{D}_{\mathrm{val}}$
			
			\STATE Initialize $\Theta$ and Adam
			
			\FOR{$e=1$ to $N_{\mathrm{epochs}}$}
			
			\FOR{each training mini-batch}
			
			\STATE $\mathbf{x}
			\leftarrow
			\left[
			\operatorname{Re}\{\widetilde{\chi}_{\mathrm{rec}}\},
			\operatorname{Im}\{\widetilde{\chi}_{\mathrm{rec}}\}
			\right]$
			
			\STATE $\widetilde{\chi}_{\mathrm{final}}
			\leftarrow
			\mathcal{U}_{\Theta}(\mathbf{x})$
			
			\STATE $\mathcal{L}_{\mathrm{total}}
			\leftarrow
			\mathcal{L}_{\mathrm{MSE}}
			\left(
			\widetilde{\chi}_{\mathrm{final}},
			\widetilde{\chi}_{\mathrm{true}}
			\right)$
			
			\STATE Backpropagate $\mathcal{L}_{\mathrm{total}}$; apply Adam updates following the gradient-accumulation rule
			
			\ENDFOR
			
			\STATE Evaluate the validation loss on $\mathcal{D}_{\mathrm{val}}$
			
			\IF{the validation loss reaches a new minimum}
			
			\STATE Save $\Theta$ as $\Theta^*$
			
			\ENDIF
			
			\IF{early stopping is enabled and triggered}
			
			\STATE \textbf{break}
			
			\ENDIF
			
			\ENDFOR
			
			\RETURN $\mathcal{U}_{\Theta^*}$
		\end{algorithmic}
	\end{algorithm}

	\subsection{Training Strategy under Noisy Conditions}
	
	Measurement noise perturbs the scattered-field data and propagates through far-field extraction and Fourier coefficient estimation.
	It can cause additional artifacts and contrast bias in the Fourier initial reconstruction.
	The noise levels used for training and testing can therefore affect U-Net reconstruction performance, alongside the network architecture and training parameters.
	
	When the training and testing noise levels are the same,
	let both noise levels be $\delta$,
	and denote the corresponding datasets by
	$\mathcal{D}_{\mathrm{train}}^{(\delta)}$ and
	$\mathcal{D}_{\mathrm{test}}^{(\delta)}$, respectively.
	The network is trained on Fourier initial reconstructions obtained from scattering data with noise level $\delta$
	and evaluated at the same noise level.
	This setting assesses reconstruction performance at a known noise level and serves as a reference for comparing training strategies.
	
	When the noise levels differ,
	a network trained at $\delta_{\mathrm{train}}$ is applied directly to testing data at
	$\delta_{\mathrm{test}}$, where
	$\delta_{\mathrm{train}}\neq\delta_{\mathrm{test}}$.
	The resulting Fourier initial reconstructions may differ from those encountered during training.
	The network may therefore yield larger reconstruction errors or additional artifacts.
	Its response to noise mismatch is assessed by fixing the training noise level and varying the testing noise level.
	The numerical experiments compare training and testing at the same noise level, direct testing under noise mismatch,
	and high-noise pretraining followed by target-noise fine-tuning.

	\section{Numerical Experiments}
	Four numerical experiments are conducted to evaluate the proposed Fourier--U-Net reconstruction method.
	They examine parameter selection, noise adaptability, complex target reconstruction, and transfer to unseen character categories.
	Experiment 1 examines the effects of the receiver distance $R$ and Fourier truncation order $N$ when training and testing use the same noise level.
	The selected configuration is used in the subsequent experiments.
	A zero-frequency compensation ablation is also performed under this configuration to assess the contribution of the {\mgq zero-frequency} term.
	Experiment 2 examines direct testing under noise mismatch and high-noise pretraining followed by target-noise fine-tuning.
	Experiment 3 evaluates reconstruction performance for single-digit and four-digit combination targets.
	Experiment 4 applies the network trained on the single-digit task in Experiment 3 to unseen character categories without retraining.
	
	In all four experiments, the square imaging domain $\Omega$ has side length $a=1$, with
	$\Omega=[-0.5,0.5]\times[-0.5,0.5].$
	The domain is discretized on a uniform $64\times64$ grid.
	The true target contrasts lie in $[0,2.5]$, and the normalization scale is therefore set to $\chi_{\mathrm{s}}=2.5$.
	
	Measurement noise is simulated by adding a complex random perturbation with a fixed relative $L_2$ norm to the scattering data.
	Let $u^s$ denote the vector of finite-distance scattered-field data corresponding to all nonzero Fourier modes.
	A complex random vector of the same dimension is generated as
	\begin{equation}
		z=\frac{\zeta_r+i\zeta_i}{\sqrt{2}},
	\end{equation}
	where all components of $\zeta_r$ and $\zeta_i$ are mutually independent and follow the standard normal distribution $\mathcal{N}(0,1)$.
	For a prescribed noise level $\delta$, the perturbed scattering data are defined by
	\begin{equation}
		u^{s,\delta}=u^s+\delta\|u^s\|_2\frac{z}{\|z\|_2}.
		\label{eq:noise_nonzero}
	\end{equation}
	When zero-frequency compensation is used, an independent perturbation at the same relative noise level is added to the theoretical far-field value $E_0$.
	This value is perturbed separately because it is not extracted from finite-distance scattered fields.
	A complex random scalar is generated as
	\begin{equation}
		z_0=\frac{\zeta_{0,r}+i\zeta_{0,i}}{\sqrt{2}},
	\end{equation}
	where $\zeta_{0,r}$ and $\zeta_{0,i}$ are independent standard normal random variables following $\mathcal{N}(0,1)$ and are independent of $z$.
	The perturbed far-field value is defined by
	\begin{equation}
		E_0^\delta=E_0+\delta|E_0|\frac{z_0}{|z_0|}.
		\label{eq:noise_zero}
	\end{equation}
	In the zero-frequency compensation ablation in Experiment 1, this compensation term is omitted as specified in the corresponding experimental setting.
	
	{\mgq In Experiments 1 and 2, each training, validation, and testing sample contains 1–3 disks.}
	Let $D_j$ denote the {\mgq $j$-th} disk, with radius $r_j$, contrast $c_j$, and center $\mathbf{y}_j$.
	The radii and contrasts are sampled from $r_j\sim\mathcal{U}(0.15,0.25)$ and $c_j\sim\mathcal{U}(0.5,2.5)$, subject to $|c_i-c_j|\geq0.3$ for $i\neq j$.
	For targets containing two or three disks, each pair is independently assigned to be separated or overlapping, with probability $0.5$ for either case.
	The distance between the centers of $D_i$ and $D_j$ is
	\begin{equation}
		d_{ij}=\|\mathbf{y}_i-\mathbf{y}_j\|_2.
	\end{equation}
	For a separated pair, the normalized separation distance is defined by
	\begin{equation}
		\beta_{ij}=\frac{d_{ij}-r_i-r_j}{\min(r_i,r_j)},
	\end{equation}
	with
	\begin{equation}
		\beta_{ij}\sim\mathcal{U}(0.2,0.4).
	\end{equation}
	For an overlapping pair, the normalized overlap is defined by
	\begin{equation}
		\alpha_{ij}=\frac{r_i+r_j-d_{ij}}{\min(r_i,r_j)},
	\end{equation}
	with
	\begin{equation}
		\alpha_{ij}\sim\mathcal{U}(0.2,0.4).
	\end{equation}
	The separation gap or overlap length is thus $20\%$--$40\%$ of the smaller disk radius.
	The radii and positions are resampled until all disks lie within $\Omega$ and satisfy the prescribed geometric constraints.
	In overlapping regions, the largest disk contrast is used to keep the target contrast within $[0,2.5]$.
	The contrast is therefore defined by
	\begin{equation}
		\chi(\mathbf{x})=
		\begin{cases}
			\displaystyle
			\max_{j:\,\mathbf{x}\in D_j}c_j,
			&
			\mathbf{x}\in\displaystyle\bigcup_jD_j,
			\\[6pt]
			0,
			&
			\mathbf{x}\notin\displaystyle\bigcup_jD_j.
		\end{cases}
		\label{eq:exp1_target_definition}
	\end{equation}
	
	Reconstruction performance is evaluated using the normalized mean squared error (NMSE) in Experiments 1 and 2.
	Experiments 3 and 4 also use the structural similarity index measure (SSIM).
	NMSE is defined as
	\begin{equation}
		\mathrm{NMSE}=\frac{\sum_{i,j}\left|\chi_{i,j}^{\mathrm{est}}-\chi_{i,j}^{\mathrm{true}}\right|^2}{\sum_{i,j}\left|\chi_{i,j}^{\mathrm{true}}\right|^2},
		\label{eq:nmse}
	\end{equation}
	where $\chi_{i,j}^{\mathrm{est}}$ denotes the direct Fourier reconstruction or the Fourier--U-Net reconstruction.
	A smaller NMSE indicates a lower reconstruction error.
	For direct Fourier reconstruction, the real part of the complex Fourier initial reconstruction is used as $\chi_{i,j}^{\mathrm{est}}$.
	No clipping is applied before computing NMSE and SSIM.
	The local SSIM is defined~\cite{wang2004ssim} as
	\begin{equation}
		\mathrm{SSIM}(x,y)=\frac{(2\mu_x\mu_y+C_1)(2\sigma_{xy}+C_2)}{(\mu_x^2+\mu_y^2+C_1)(\sigma_x^2+\sigma_y^2+C_2)},
		\label{eq:ssim}
	\end{equation}
	where $x$ and $y$ denote corresponding local patches of the reconstructed and true target images.
	The quantities $\mu_x,\mu_y$ are the local means, $\sigma_x,\sigma_y$ are the local standard deviations, and $\sigma_{xy}$ is the local covariance.
	The constants are $C_1=(K_1L)^2$ and $C_2=(K_2L)^2$, with $K_1=0.01$, $K_2=0.03$, and a fixed contrast range $L=2.5$.
	The local SSIM scores are averaged to obtain the image-level value.
	A value closer to $1$ indicates higher structural similarity.
	Unless otherwise specified, the main result tables report the mean and sample standard deviation across testing samples.
	For Fourier--U-Net, evaluation metrics are first computed separately for each model and testing sample.
	For each sample, the metric values are then averaged across models.
	The reported mean and sample standard deviation are calculated across these per-sample metric averages.
	For direct Fourier reconstruction, the mean and sample standard deviation are calculated over the same testing samples.
	The statistics for parameter selection in Experiment 1 are specified separately.
	
	The experiments are implemented using PyTorch 2.12.1 and Python 3.11.12.
	Training and inference are performed on an NVIDIA GeForce RTX 5070 Laptop GPU with $8$ GB of memory.

	\subsection{Experiment 1: Fourier Truncation Order Selection and Reconstruction Performance Comparison}
		Experiment 1 first selects the receiver radius $R$ for finite-distance measurements and then determines the Fourier truncation order $N$ at this radius.
	For nonzero Fourier modes, the relative far-field error is defined as
	\begin{equation}
		e_{\mathrm{ff}}=\frac{\left\|E_{\infty}^{(R)}-E_{\infty}^{\mathrm{ref}}\right\|_2}{\left\|E_{\infty}^{\mathrm{ref}}\right\|_2},
		\label{eq:exp1_farfield_error}
	\end{equation}
	where $E_{\infty}^{(R)}$ is the far-field pattern extracted from finite-distance scattered fields using the asymptotic relation.
	The reference pattern $E_{\infty}^{\mathrm{ref}}$ is computed directly from the Born model.
	The zero-frequency contribution is supplied independently by the theoretical far-field value at the small wave number $k_0=0.01$.
	Since it is not extracted from finite-distance data, it is excluded from the receiver-radius error assessment.
	
	The candidate receiver radii are $R=3,5,10,20$.
	A total of $300$ targets are generated, comprising $100$ single-disk, $100$ double-disk, and $100$ triple-disk targets.
	For each receiver radius and truncation order, the relative far-field error is calculated for every target.
	Table~\ref{tab:exp1_farfield_validation} reports the mean and standard deviation of these errors over the $300$ targets.
	\begin{table}[htbp]
		\centering
		\caption{Far-field relative errors at different receiver radii}
		\label{tab:exp1_farfield_validation}
		\renewcommand{\arraystretch}{1.15}
		\setlength{\tabcolsep}{8pt}
		\fontsize{9.5pt}{11.5pt}\selectfont
		\begin{tabular}{c|cccc}
			\toprule
			$N$ & $R=3$ & $R=5$ & $R=10$ & $R=20$ \\
			\midrule
			3
			& $0.0662\pm0.0255$
			& $0.0397\pm0.0153$
			& $0.0198\pm0.0077$
			& $0.0099\pm0.0038$ \\
			4
			& $0.0855\pm0.0366$
			& $0.0513\pm0.0220$
			& $0.0256\pm0.0110$
			& $0.0128\pm0.0055$ \\
			5
			& $0.1032\pm0.0449$
			& $0.0619\pm0.0270$
			& $0.0310\pm0.0135$
			& $0.0155\pm0.0068$ \\
			10
			& $0.2013\pm0.0854$
			& $0.1216\pm0.0521$
			& $0.0610\pm0.0262$
			& $0.0305\pm0.0131$ \\
			\bottomrule
		\end{tabular}
	\end{table}
	
	Let $\overline{e}_{\mathrm{ff}}(N,R)$ denote the mean relative far-field error over the $300$ targets at truncation order $N$ and receiver radius $R$.
	Receiver-radius selection is based on $N=3,4,5$, the candidate orders for the Fourier--U-Net configuration.
	The mean error across these three orders is defined by
	\begin{equation}
		\overline{e}_{\mathrm{ff}}(R)=\frac{1}{3}\sum_{N\in\{3,4,5\}}\overline{e}_{\mathrm{ff}}(N,R),
	\end{equation}
	with the requirement
	\begin{equation}
		\overline{e}_{\mathrm{ff}}(R)\leq0.1.
	\end{equation}
	According to Table~\ref{tab:exp1_farfield_validation}, $R=3$ is the smallest candidate radius satisfying this criterion.
	It is therefore used in the subsequent experiments.
	
	With $R=3$ fixed, the candidate truncation orders $N=3,4,5$ are compared using reconstruction performance.
	The training and validation sets contain $3000$ and $400$ samples, respectively.
	An independent parameter selection set contains $100$ single-disk, $100$ double-disk, and $100$ triple-disk targets, giving $300$ samples in total.
	The same target samples are used across truncation orders to control for differences in target geometry.
	For each order, three U-Nets are trained independently with different random seeds to reduce sensitivity to network initialization.
	
	Training and parameter selection use a relative noise level of $5\%$.
	Noise generation, network architecture, normalization, loss function, and optimizer follow the settings described above.
	The physical batch size is $4$, and accumulating gradients over $8$ mini-batches gives an effective batch size of $32$.
	The initial learning rate is $10^{-4}$ and is halved every $20$ epochs.
	Training runs for at most $100$ epochs, with an early-stopping patience of $15$ epochs.
	The parameters with the lowest validation MSE are retained.
	For each random seed, NMSE is first averaged over the $300$ parameter selection samples.
	Table~\ref{tab:exp1_selection} reports the mean and standard deviation of these three model-level mean NMSE values.
	
	\begin{table}[htbp]
		\centering
		\caption{Selection results for Fourier truncation orders at $5\%$ noise}
		\label{tab:exp1_selection}
		\renewcommand{\arraystretch}{1.15}
		\setlength{\tabcolsep}{8pt}
		\fontsize{9.5pt}{11.5pt}\selectfont
		\begin{tabular}{c|ccc}
			\toprule
			Order $N$ & 3 & 4 & 5 \\
			\midrule
			Modes $M=(2N+1)^2$
			& 49 & 81 & 121 \\
			Selection NMSE
			& $0.0132\pm0.0009$
			& $0.0104\pm0.0013$
			& $0.0075\pm0.0004$ \\
			Mean far-field error at $R=3$
			& 0.0662 & 0.0855 & 0.1032 \\
			\bottomrule
		\end{tabular}
	\end{table}
	
	Across the tested truncation orders, the selection NMSE decreases as $N$ increases, while the finite-distance far-field extraction error increases.
	For truncation-order selection, candidates must satisfy
	$\overline{e}_{\mathrm{ff}}(N,3)\leq0.1$.
	Among the eligible candidates, $N=4$ gives the lowest selection NMSE.
	The subsequent experiments therefore use $R=3$ and $N=4$.
	
\begin{figure}[htbp]
		\centering
		\includegraphics[
		width=1.0\textwidth,
		keepaspectratio
		]{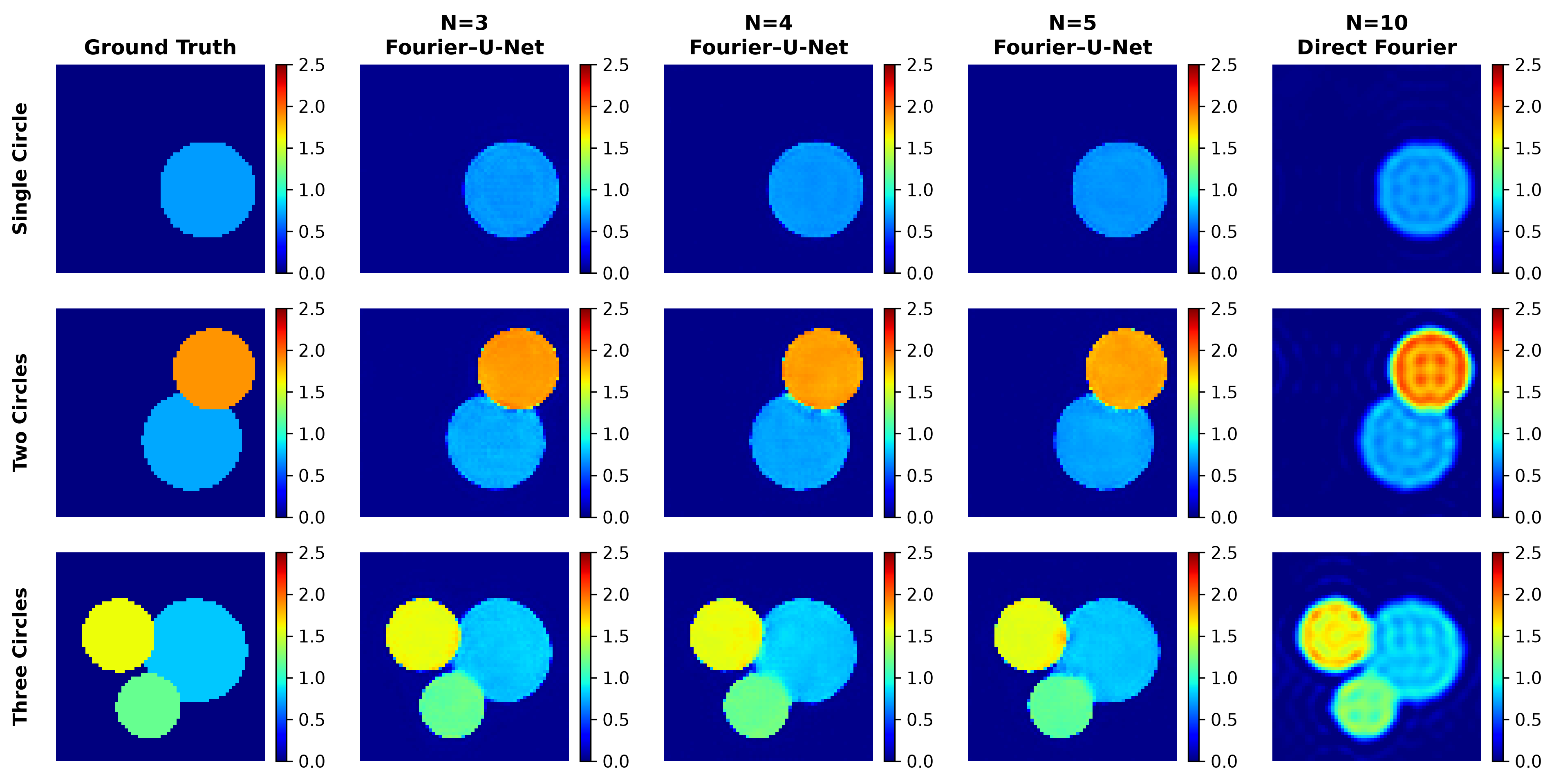}
		\caption{Reconstruction results for different Fourier truncation orders on the independent final testing set under $5\%$ noise. The three rows correspond to single-disk, double-disk, and triple-disk targets, respectively; from left to right are the ground truth,
			$N=3$ Fourier--U-Net, $N=4$ Fourier--U-Net, $N=5$ Fourier--U-Net, and $N=10$ direct Fourier reconstruction.}
		\label{fig:exp1_generalization}
	\end{figure}

	After selecting $R$ and $N$, a final testing set is generated independently of the training, validation, and parameter selection sets.
	It contains $100$ single-disk, $100$ double-disk, and $100$ triple-disk targets, giving $300$ samples in total.
	The testing noise level is $5\%$.
	This set is not used to select $R$, $N$, or the training parameters.
	Figure~\ref{fig:exp1_generalization} shows the Fourier--U-Net reconstructions for $N=3,4,5$ on this set.
	Direct Fourier reconstruction with $N=10$ is included as a baseline using more Fourier modes.
	
	For the targets shown in Fig.~\ref{fig:exp1_generalization}, Fourier--U-Net with $N=3,4,5$ recovers the number, locations, and main shapes of the disks.
	The $N=10$ direct Fourier reconstructions exhibit oscillations and local artifacts despite retaining more Fourier modes.
	Table~\ref{tab:exp1_unet_fourier_nmse} compares NMSE for $N=4$ direct Fourier reconstruction, $N=4$ Fourier--U-Net, and $N=10$ direct Fourier reconstruction across the three target types.
	\begin{table}[htbp]
		\centering
		\caption{NMSE on the independent final test set at $5\%$ noise}
		\label{tab:exp1_unet_fourier_nmse}
		\renewcommand{\arraystretch}{1.15}
		\setlength{\tabcolsep}{9pt}
		\fontsize{9.5pt}{11.5pt}\selectfont
		\begin{tabular}{c|ccc}
			\toprule
			Target type
			& $N=4$ direct Fourier
			& $N=4$ Fourier--U-Net
			& $N=10$ direct Fourier \\
			\midrule
			Single disk
			& $0.1052\pm0.0207$
			& $0.0100\pm0.0073$
			& $0.0445\pm0.0072$ \\
			Double disk
			& $0.0988\pm0.0156$
			& $0.0089\pm0.0030$
			& $0.0414\pm0.0056$ \\
			Triple disk
			& $0.0955\pm0.0170$
			& $0.0106\pm0.0025$
			& $0.0394\pm0.0058$ \\
			\bottomrule
		\end{tabular}
	\end{table}
	
	For all three target types in Table~\ref{tab:exp1_unet_fourier_nmse}, $N=4$ Fourier--U-Net has lower NMSE than direct Fourier reconstruction with either $N=4$ or $N=10$.
	Under the Born approximation, with finite-distance measurements and $5\%$ noise used here, U-Net refinement improves the initial reconstruction at $N=4$.
	The comparison between $N=4$ Fourier--U-Net and $N=10$ direct Fourier reconstruction does not isolate the effect of truncation order.
	
	A zero-frequency compensation ablation is performed at receiver radius $R=3$, truncation order $N=4$, and noise level $5\%$.
	With compensation, the zero-frequency coefficient is approximated using the theoretical far-field value at a small wave number.
	Without compensation, the center entry of the Fourier coefficient matrix is set to $\widehat{\chi}_{\boldsymbol{0}}=0$.
	The nonzero Fourier modes, noise realizations, network architecture, and training parameters remain unchanged.
	
	The ablation uses the same independent final testing set of $300$ targets, with $100$ targets of each disk configuration.
	For the case without compensation, three $N=4$ U-Nets are retrained using different random seeds.
	Direct Fourier reconstruction and Fourier--U-Net are evaluated on the same testing targets.
	For Fourier--U-Net, NMSE is computed separately for each model and testing target.
	For each target, the three NMSE values are averaged across models.
	Table~\ref{tab:zero_frequency_ablation} reports the mean and sample standard deviation of these per-target NMSE averages for Fourier--U-Net.
	For direct Fourier reconstruction, these statistics are computed from the NMSE values of the individual targets.
	\begin{table}[htbp]
		\centering
		\caption{Zero-frequency ablation on the independent test set ($R=3$, $N=4$, $5\%$ noise)}
		\label{tab:zero_frequency_ablation}
		\renewcommand{\arraystretch}{1.15}
		\setlength{\tabcolsep}{6pt}
		\fontsize{9.5pt}{11.5pt}\selectfont
		\begin{tabular}{@{}lccc@{}}
			\toprule
			Method
			& \shortstack{Zero-frequency compensation}
			& \shortstack{Available coefficients}
			& NMSE \\
			\midrule
			Direct Fourier
			& With
			& 81
			& $0.0998\pm0.0183$ \\
			Direct Fourier
			& Without
			& 80
			& $0.3218\pm0.0820$ \\
			Fourier--U-Net
			& With
			& 81
			& $0.0098\pm0.0048$ \\
			Fourier--U-Net
			& Without
			& 80
			& $0.0105\pm0.0053$ \\
			\bottomrule
		\end{tabular}
	\end{table}
	
	As shown in Table~\ref{tab:zero_frequency_ablation}, removing zero-frequency compensation increases the NMSE of direct Fourier reconstruction from $0.0998$ to $0.3218$.
	The compensation term therefore has a substantial effect on the direct Fourier reconstructions of the tested disk targets.
	
	For Fourier--U-Net, NMSE increases from $0.0098$ to $0.0105$ when compensation is removed and the network is retrained.
	The latter remains lower than the NMSE of direct Fourier reconstruction without compensation.
	Under the Born approximation, with finite-distance measurements and $5\%$ noise used here, this improvement is retained when training and testing both omit the compensation term.

	\subsection{Experiment 2: High-Noise Pretraining and Target-Noise Fine-Tuning}
		With the receiver radius $R=3$ and Fourier truncation order $N=4$ fixed, Experiment 2 examines Fourier--U-Net reconstruction when training and testing noise levels differ.
	It also compares high-noise pretraining followed by target-noise fine-tuning with training from scratch.
	
\begin{figure}[htbp]
		\centering
		\includegraphics[
		width=0.6\textwidth,
		keepaspectratio
		]{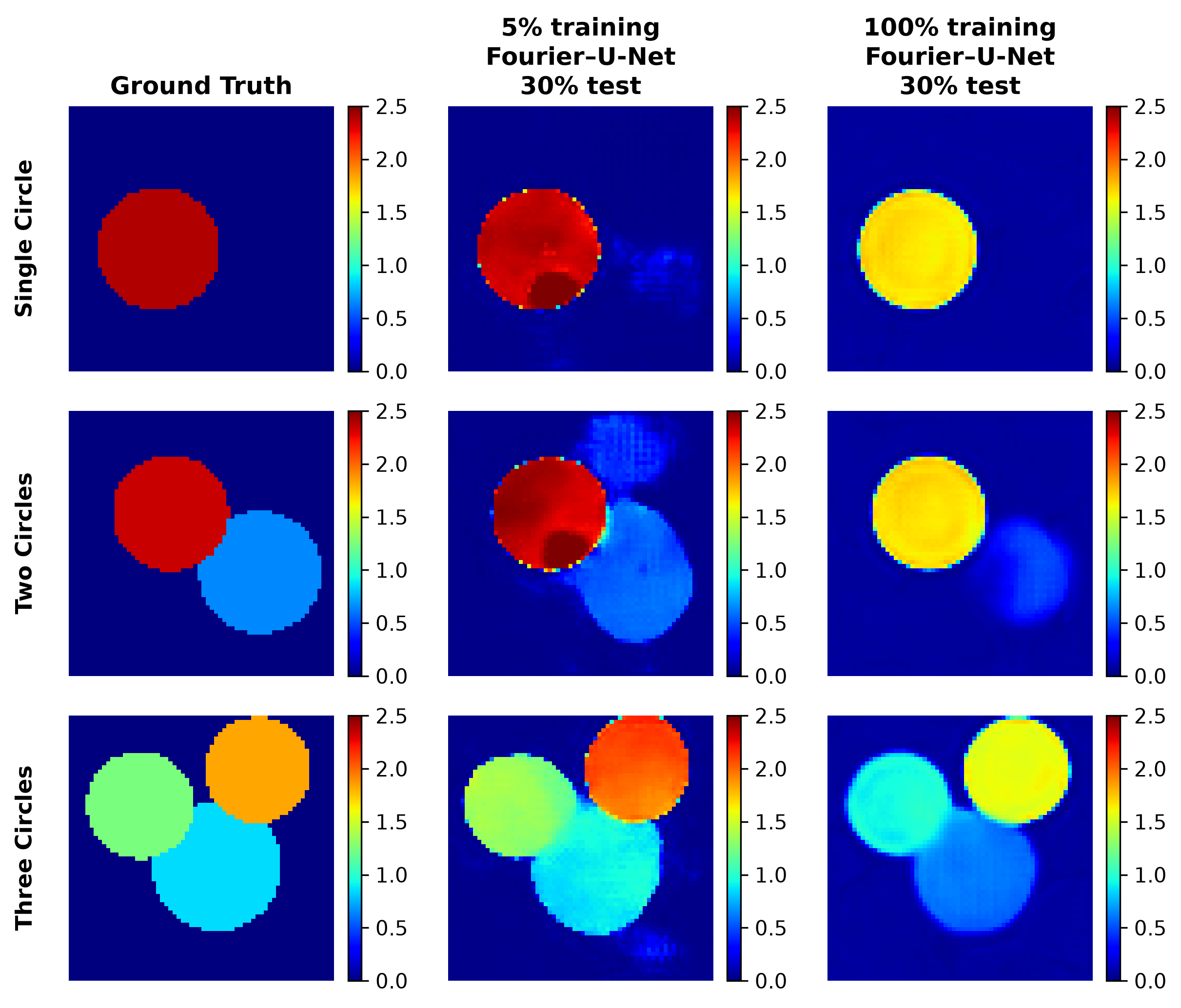}
		\caption{Reconstruction results of models trained at $5\%$ and $100\%$ noise, respectively, and directly applied to test data at $30\%$ noise ($R=3$, $N=4$).}
		\label{fig:exp2_direct30}
	\end{figure}

	The data generation method, measurement procedure, and network architecture follow Experiment 1.
	Two models, trained at $5\%$ and $100\%$ noise, respectively, are first applied directly to testing data at $30\%$ noise without fine-tuning.
	Both models are evaluated on the same single-disk, double-disk, and triple-disk targets with identical noise realizations.
	The reconstructions are shown in Fig.~\ref{fig:exp2_direct30}.
	
	For the examples in Fig.~\ref{fig:exp2_direct30}, the model trained at $5\%$ noise produces spurious structures and local artifacts near the targets.
	The model trained at $100\%$ noise recovers the number of disks and their main contours, but shows contrast bias within the targets.
	These examples illustrate the artifacts and contrast errors that can occur under noise mismatch.
	
	Three training strategies are compared: Transfer, Scratch, and Direct100.
	Transfer begins with pretraining for $60$ epochs at $100\%$ noise.
	The checkpoint with the lowest validation MSE initializes separate fine-tuning runs of $40$ epochs at target noise levels of $5\%$, $30\%$, and $50\%$.
	The initial learning rate is $1\times10^{-4}$ for pretraining and $3\times10^{-5}$ for fine-tuning.
	It is halved every $20$ epochs during pretraining and every $10$ epochs during fine-tuning.
	Scratch is trained from random initialization for $100$ epochs at each target noise level.
	Its initial learning rate is $1\times10^{-4}$ and is halved every $20$ epochs.
	Direct100 is trained for $100$ epochs only at $100\%$ noise and tested directly at the target noise levels without fine-tuning.
	
	For a single target noise task, Transfer uses $60$ pretraining epochs and $40$ fine-tuning epochs, matching the $100$ epochs used by Scratch and Direct100.
	All training stages run for the prescribed fixed numbers of epochs, without early stopping.
	The checkpoint with the lowest validation MSE is retained at each stage.
	Each strategy is evaluated using three different random seeds.
	For each seed, the three Transfer fine-tuning runs start from the same pretrained checkpoint.
	
	An independent testing set contains $100$ single-disk, $100$ double-disk, and $100$ triple-disk targets, giving $300$ samples in total.
	At each noise level, all three strategies are evaluated on the same testing targets and noise realizations.
	This set is not used for training, validation, or parameter updates.
	Table~\ref{tab:exp2_transfer_nmse} reports the NMSE results at the three target noise levels.
	The relative NMSE reductions of Transfer with respect to Scratch and Direct100 are denoted by $I_{\mathrm{T/S}}$ and $I_{\mathrm{T/D}}$, respectively.
	These rates are calculated from the unrounded mean NMSE values.
	
	The training time comparison covers the three target noise tasks and all three random seeds.
	Transfer uses $3$ shared pretraining runs and $9$ fine-tuning runs, while Scratch requires $9$ independent training runs.
	Table~\ref{tab:exp2_multinoise_runtime} reports the stage times and total times for the two strategies.
	The total times are obtained by summing the displayed stage times.
	
	\begin{table}[htbp]
		\centering
		\caption{NMSE of three training strategies at different noise levels}
		\label{tab:exp2_transfer_nmse}
		\renewcommand{\arraystretch}{1.15}
		\setlength{\tabcolsep}{10.5pt}
		\fontsize{9.5pt}{11.5pt}\selectfont
		\begin{tabular}{@{}cccccc@{}}
			\toprule
			Noise
			& Scratch
			& Transfer
			& Direct100
			& $I_{\mathrm{T/S}}$
			& $I_{\mathrm{T/D}}$ \\
			\midrule
			$5\%$
			& $0.0095\pm0.0047$
			& $0.0105\pm0.0056$
			& $0.1043\pm0.0245$
			& $-10.61\%$
			& $89.98\%$ \\
			$30\%$
			& $0.0258\pm0.0098$
			& $0.0255\pm0.0109$
			& $0.0971\pm0.0259$
			& $1.21\%$
			& $73.78\%$ \\
			$50\%$
			& $0.0410\pm0.0160$
			& $0.0394\pm0.0157$
			& $0.0881\pm0.0277$
			& $3.96\%$
			& $55.25\%$ \\
			\bottomrule
		\end{tabular}
	\end{table}
	
	\begin{table}[htbp]
		\centering
		\caption{Training time comparison of Transfer and Scratch for three target noise levels}
		\label{tab:exp2_multinoise_runtime}
		\renewcommand{\arraystretch}{1.15}
		\setlength{\tabcolsep}{10pt}
		\fontsize{9.5pt}{11.5pt}\selectfont
		\begin{tabular}{@{}lcccc@{}}
			\toprule
			Strategy
			& Pretraining ($\mathrm{s}$)
			& Target-noise training ($\mathrm{s}$)
			& Total time ($\mathrm{s}$)
			& Time reduction \\
			\midrule
			Transfer
			& 1811.6
			& 3634.0
			& 5445.6
			& $39.08\%$ \\
			Scratch
			& --
			& 8938.8
			& 8938.8
			& -- \\
			\bottomrule
		\end{tabular}
	\end{table}
	
	As shown in Table~\ref{tab:exp2_transfer_nmse}, Transfer has a higher mean NMSE than Scratch at $5\%$ noise, with values of $0.0105$ and $0.0095$, respectively.
	At $30\%$ and $50\%$ noise, Transfer gives numerically slightly lower NMSE values of $0.0255$ and $0.0394$, compared with $0.0258$ and $0.0410$ for Scratch.
	For synthetic scattering data generated under the Born approximation, with the fixed relative noise model used here, high-noise pretraining followed by target-noise fine-tuning achieves reconstruction accuracy close to Scratch.
	Direct100 has higher mean NMSE at all three noise levels.
	Relative to Direct100, Transfer reduces NMSE by $89.98\%$, $73.78\%$, and $55.25\%$, respectively.
	Target-noise fine-tuning therefore improves reconstruction compared with direct testing of the model trained only at $100\%$ noise.
	
\begin{figure}[htbp]
		\centering
		\includegraphics[
		width=0.8\textwidth,
		keepaspectratio
		]{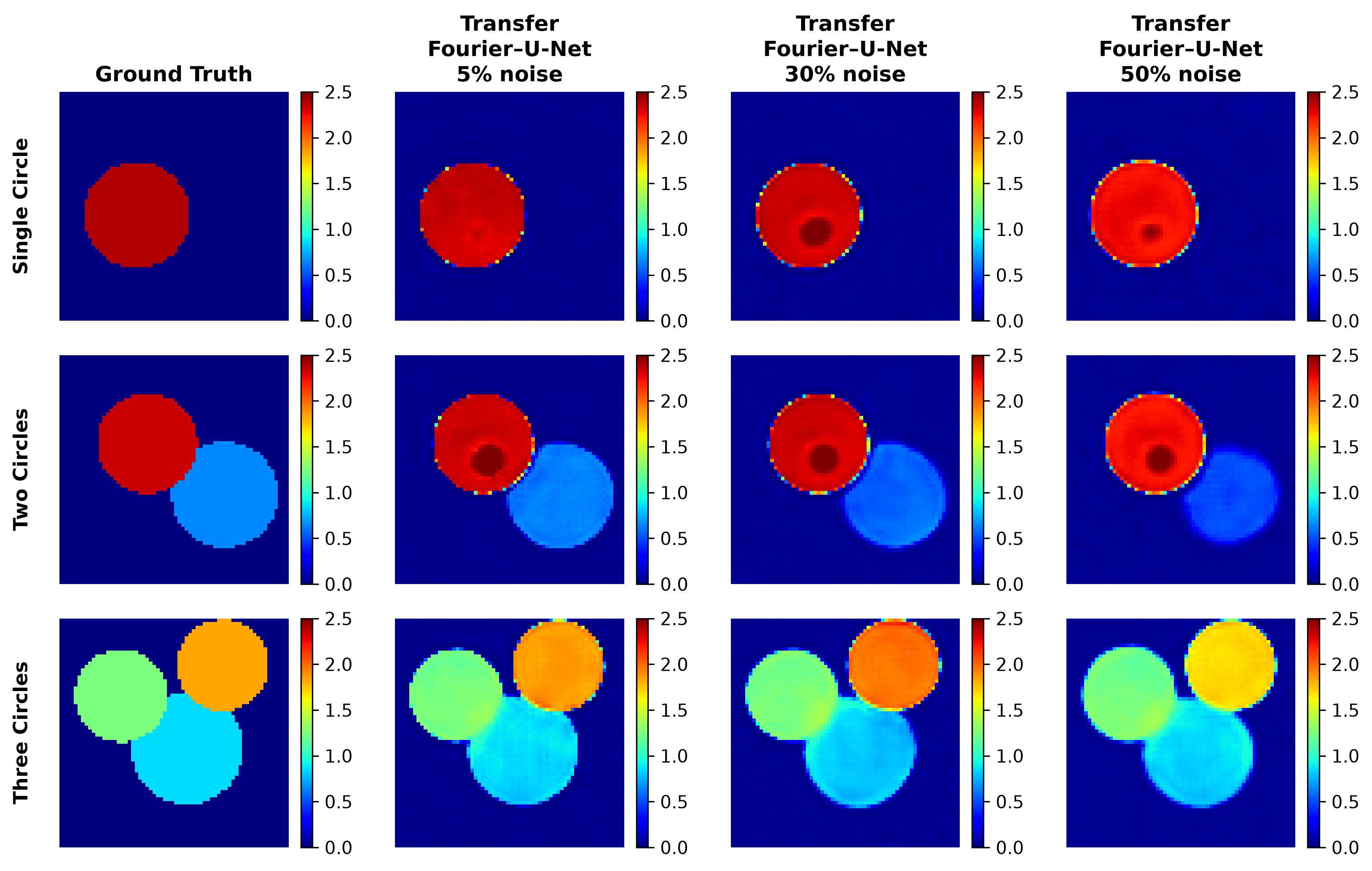}
		\caption{Reconstruction results after high-noise pretraining and target-noise fine-tuning at different noise levels ($R=3$, $N=4$).}
		\label{fig:exp2_transfer}
	\end{figure}

	Across the three target noise tasks, Transfer takes $5445.6\,\mathrm{s}$, compared with $8938.8\,\mathrm{s}$ for Scratch, as reported in Table~\ref{tab:exp2_multinoise_runtime}.
	This corresponds to a saving of $3493.2\,\mathrm{s}$, or $39.08\%$.
	For each random seed, Transfer uses one $60$-epoch pretraining run and three $40$-epoch fine-tuning runs, while Scratch uses three independent $100$-epoch training runs.
	The time reduction reflects this shared pretraining schedule and its smaller total number of training epochs.
	Under these schedules, Transfer reduces total training time while attaining mean NMSE values close to those of Scratch.
	Figure~\ref{fig:exp2_transfer} shows the Transfer reconstructions after fine-tuning at target noise levels of $5\%$, $30\%$, and $50\%$.
	
	For the displayed targets, Transfer recovers the number, locations, and main shapes of the disks at all three noise levels.
	As the noise level increases from $5\%$ to $50\%$, the boundaries become smoother and the contrast errors within the disks increase.
	The main separation and overlap relationships in the double-disk and triple-disk targets are preserved.

	\subsection{Experiment 3: Reconstruction of Complex-Structure Targets}
		Experiment 3 evaluates Fourier--U-Net on two-dimensional contrast fields constructed from MNIST handwritten digits~\cite{lecun1998gradient}.
	The measurements use $R=3$ and $N=4$, as selected in Experiment 1, with a relative noise level of $5\%$.
	For a normalized grayscale image $I(\mathbf{x})\in[0,1]$, the target contrast is defined as $\chi(\mathbf{x})=2.5I(\mathbf{x})$.
	The contrast upper bound is therefore $2.5$, as in the preceding disk experiments.
	
	Two tasks are considered: single digits and four-digit combinations.
	For single digits, the original $28\times28$ images are interpolated to $64\times64$.
	The training, validation, and testing sets contain $5000$, $1000$, and $1000$ samples, respectively.
	The testing set includes $100$ samples of each digit from $0$ to $9$.
	For four-digit combinations, each digit is resized to $32\times32$, and four digits are arranged in a $2\times2$ layout to form a $64\times64$ target.
	The corresponding training, validation, and testing sets contain $1500$, $250$, and $250$ samples, respectively.
	Separate U-Nets are trained for the two tasks, with three independent training runs using different random seeds for each task.
	The remaining training settings follow those of Experiment 1.
	The testing sets are independent of the corresponding training and validation sets and are not used for model selection.
	Table~\ref{tab:exp3_results} summarizes the reconstruction performance for both tasks.
	
	\begin{table}[htbp]
		\centering
		\caption{Reconstruction performance for MNIST complex-structure targets}
		\label{tab:exp3_results}
		\renewcommand{\arraystretch}{1.15}
		\setlength{\tabcolsep}{7pt}
		\fontsize{9.5pt}{11.5pt}\selectfont
		\begin{tabular}{lcccc}
			\toprule
			Task
			& \multicolumn{2}{c}{Direct Fourier}
			& \multicolumn{2}{c}{Fourier--U-Net} \\
			\cmidrule(lr){2-3}\cmidrule(lr){4-5}
			& NMSE & SSIM & NMSE & SSIM \\
			\midrule
			Single digit
			& $0.1448\pm0.0613$
			& $0.3810\pm0.0626$
			& $0.0119\pm0.0078$
			& $0.9216\pm0.0101$ \\
			Four-digit combination
			& $0.3926\pm0.0505$
			& $0.2601\pm0.0246$
			& $0.1213\pm0.0389$
			& $0.5225\pm0.0402$ \\
			\bottomrule
		\end{tabular}
	\end{table}
	
\begin{figure}[htbp]
		\centering
		\includegraphics[
		width=0.9\textwidth,
		keepaspectratio
		]{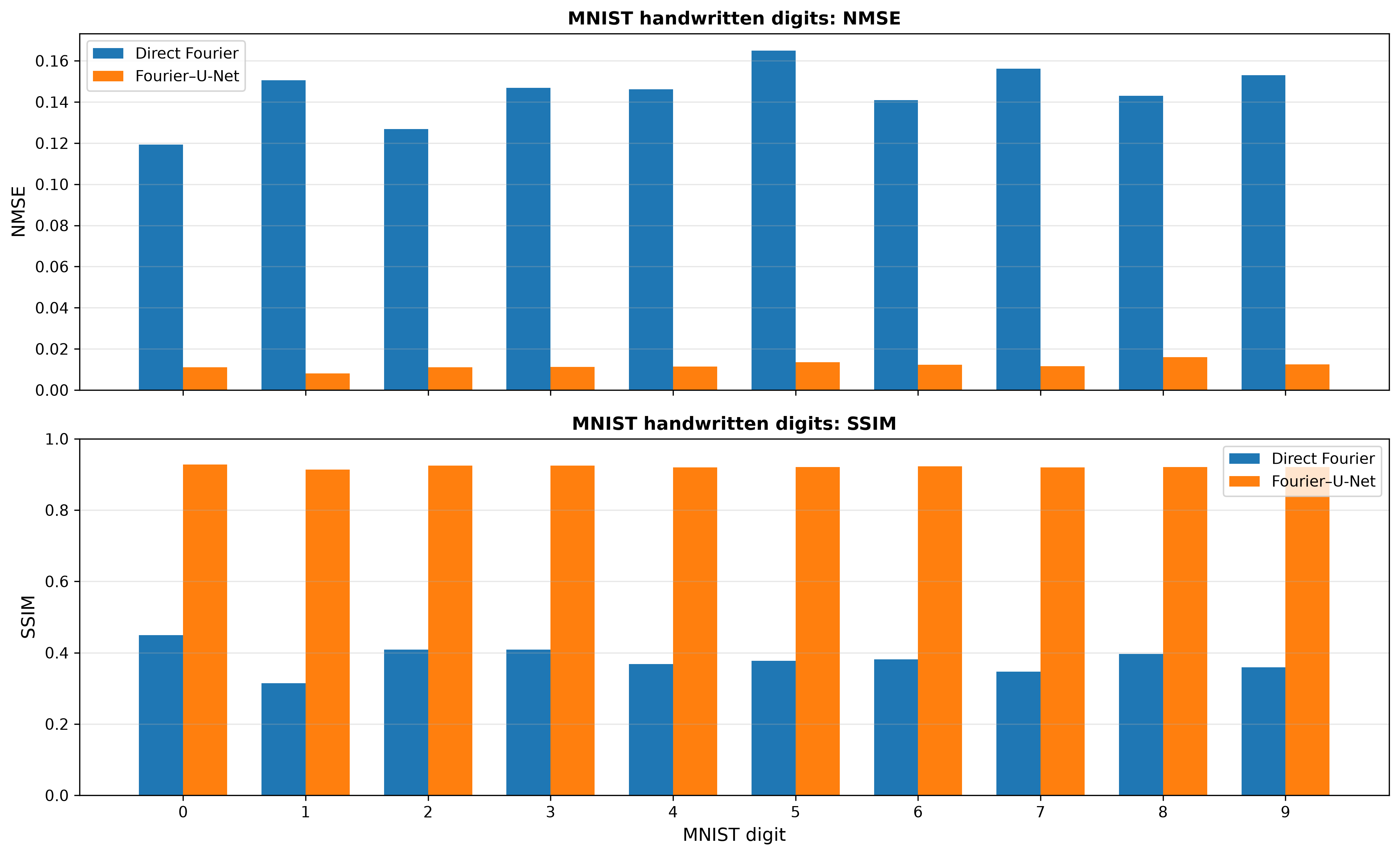}
		\caption{Comparison of NMSE and SSIM for different MNIST handwritten digits}
		\label{fig:exp3_digit_metrics}
	\end{figure}

	For single digits, Fourier--U-Net reduces mean NMSE from $0.1448$ to $0.0119$ and increases mean SSIM from $0.3810$ to $0.9216$, as shown in Table~\ref{tab:exp3_results}.
	For four-digit combinations, NMSE decreases from $0.3926$ to $0.1213$, while SSIM increases from $0.2601$ to $0.5225$.
	These results show that Fourier--U-Net improves reconstruction of the tested MNIST targets for synthetic scattering data generated under the Born approximation.
	The four-digit task has higher NMSE and lower SSIM than the single-digit task, indicating greater reconstruction difficulty under these settings.
	Figure~\ref{fig:exp3_digit_metrics} shows the mean NMSE and SSIM for each MNIST digit from $0$ to $9$.
	
\begin{figure}[htbp]
		\centering
		\includegraphics[
		width=0.95\textwidth,
		keepaspectratio
		]{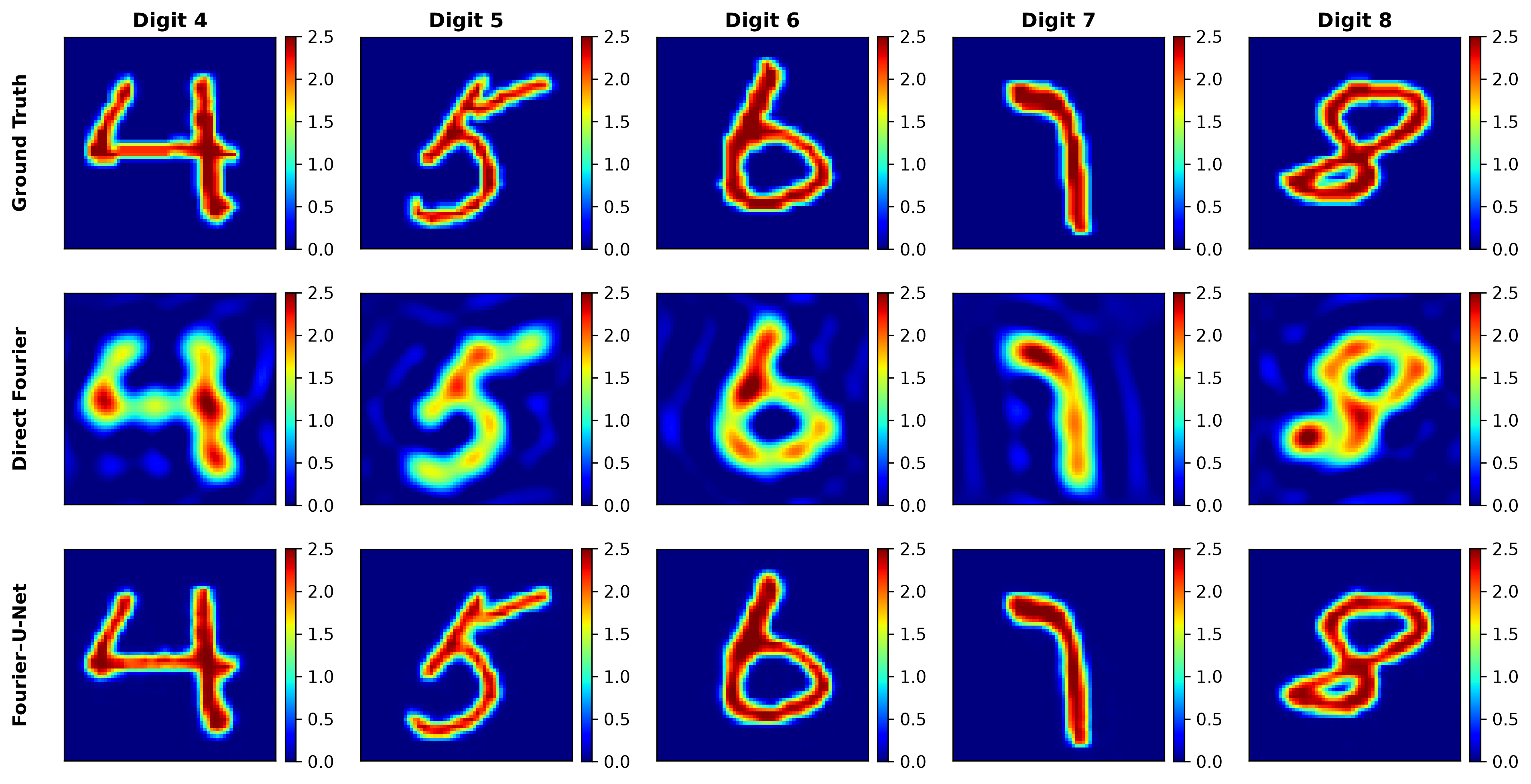}
		\caption{Reconstruction results for individual MNIST digits}
		\label{fig:exp3_single}
	\end{figure}

	For every digit category in Fig.~\ref{fig:exp3_digit_metrics}, Fourier--U-Net has lower mean NMSE and higher mean SSIM than direct Fourier reconstruction.
	The differences between categories may reflect variations in stroke shapes and internal structures.
	In Fig.~\ref{fig:exp3_single}, direct Fourier reconstruction captures the overall digit shapes but produces blurred strokes and local artifacts.
	U-Net refinement better recovers the contours, stroke widths, and internal holes, giving results closer to the true targets.
	
\begin{figure}[htbp]
		\centering
		\includegraphics[
		width=0.7\textwidth,
		keepaspectratio
		]{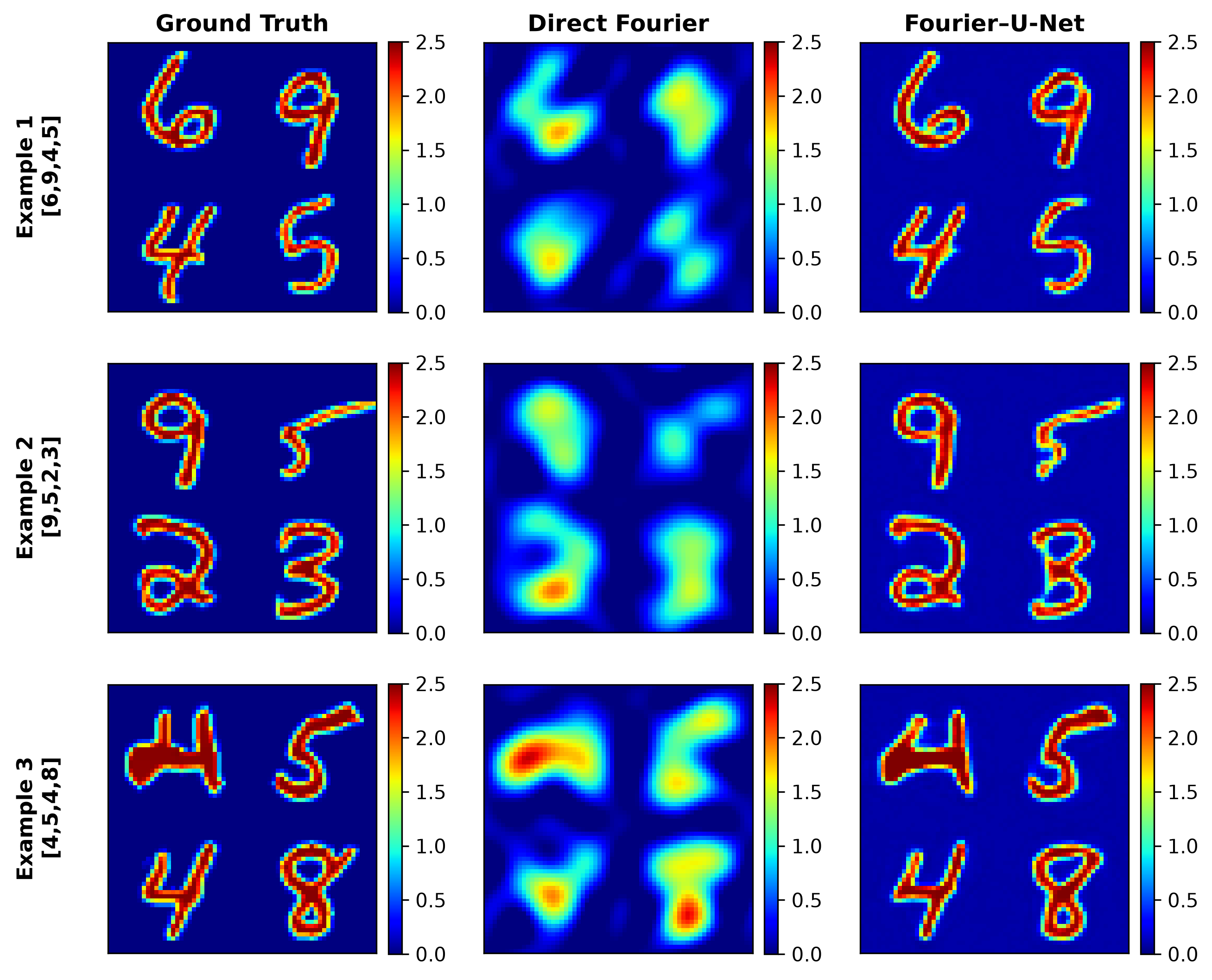}
		\caption{Reconstruction results for $2\times2$ MNIST digit combinations}
		\label{fig:exp3_combo}
	\end{figure}

	In Fig.~\ref{fig:exp3_combo}, direct Fourier reconstruction produces smoothed strokes and distorted shapes, while Fourier--U-Net recovers the digit positions and main stroke structures.
	Some local shape deviations and contrast errors remain in the four-digit reconstructions, consistent with their higher NMSE and lower SSIM.
	For synthetic scattering data generated under the Born approximation, Fourier--U-Net improves direct Fourier reconstruction of the tested irregular targets with complex boundaries and internal structures.
	
	\subsection{Experiment 4: Transfer to Unseen Character Categories}
	Experiment 4 evaluates the transfer of Fourier--U-Net to unseen character categories using the U-Nets trained on the single-digit MNIST task in Experiment 3.
	All network parameters remain frozen during testing, without retraining or fine-tuning.
	The measurement configuration, contrast normalization, and noise settings follow Experiment 3, with $R=3$, $N=4$, and a relative noise level of $5\%$.
	
	Two testing datasets are used.
	From the EMNIST handwritten letter dataset~\cite{cohen2017emnist}, $100$ samples are selected from each of $26$ categories, giving $2600$ testing samples.
	From the Chinese handwritten numeral dataset~\cite{nazarpour2017chinese}, $100$ samples are selected from each of nine categories representing one to nine, giving $900$ samples.
	EMNIST images are corrected for orientation and resized to $64\times64$.
	Chinese numeral images are converted to bright strokes on a dark background, cropped to the foreground region, resized proportionally to a longest side of $48$ pixels, and centered on a $64\times64$ grid.
	Neither dataset is used for training, validation, or parameter selection in Experiment 3.
	The same measurement and reconstruction procedures are applied using the three frozen U-Nets from that experiment.
	Table~\ref{tab:exp4_results} compares direct Fourier reconstruction at $N=4$ with Fourier--U-Net.
	The reported means and sample standard deviations follow the statistical procedure specified for the main result tables.
	
	\begin{table}[htbp]
		\centering
		\caption{Reconstruction performance for unseen character-category targets}
		\label{tab:exp4_results}
		\renewcommand{\arraystretch}{1.15}
		\setlength{\tabcolsep}{10pt}
		\fontsize{9.5pt}{11.5pt}\selectfont
		\begin{tabular}{@{}lcccc@{}}
			\toprule
			Test target
			& \multicolumn{2}{c}{Direct Fourier}
			& \multicolumn{2}{c}{Fourier--U-Net} \\
			\cmidrule(lr){2-3}\cmidrule(lr){4-5}
			& NMSE & SSIM & NMSE & SSIM \\
			\midrule
			EMNIST letters
			& $0.0907\pm0.0596$
			& $0.4748\pm0.0993$
			& $0.0164\pm0.0102$
			& $0.9079\pm0.0218$ \\
			\shortstack[l]{Chinese numerals}
			& $0.3880\pm0.1202$
			& $0.4847\pm0.0791$
			& $0.2635\pm0.1119$
			& $0.7828\pm0.0745$ \\
			\bottomrule
		\end{tabular}
	\end{table}
	
	On the EMNIST letter testing set, Fourier--U-Net reduces mean NMSE from $0.0907$ to $0.0164$ and increases mean SSIM from $0.4748$ to $0.9079$.
	Although trained only on MNIST digits, the frozen networks improve reconstruction of the tested English letters.
	This indicates transfer to unseen letter categories for synthetic scattering data generated under the Born approximation at the same contrast scale.
	
\begin{figure}[htbp]
		\centering
		\includegraphics[
		width=0.95\textwidth,
		keepaspectratio
		]{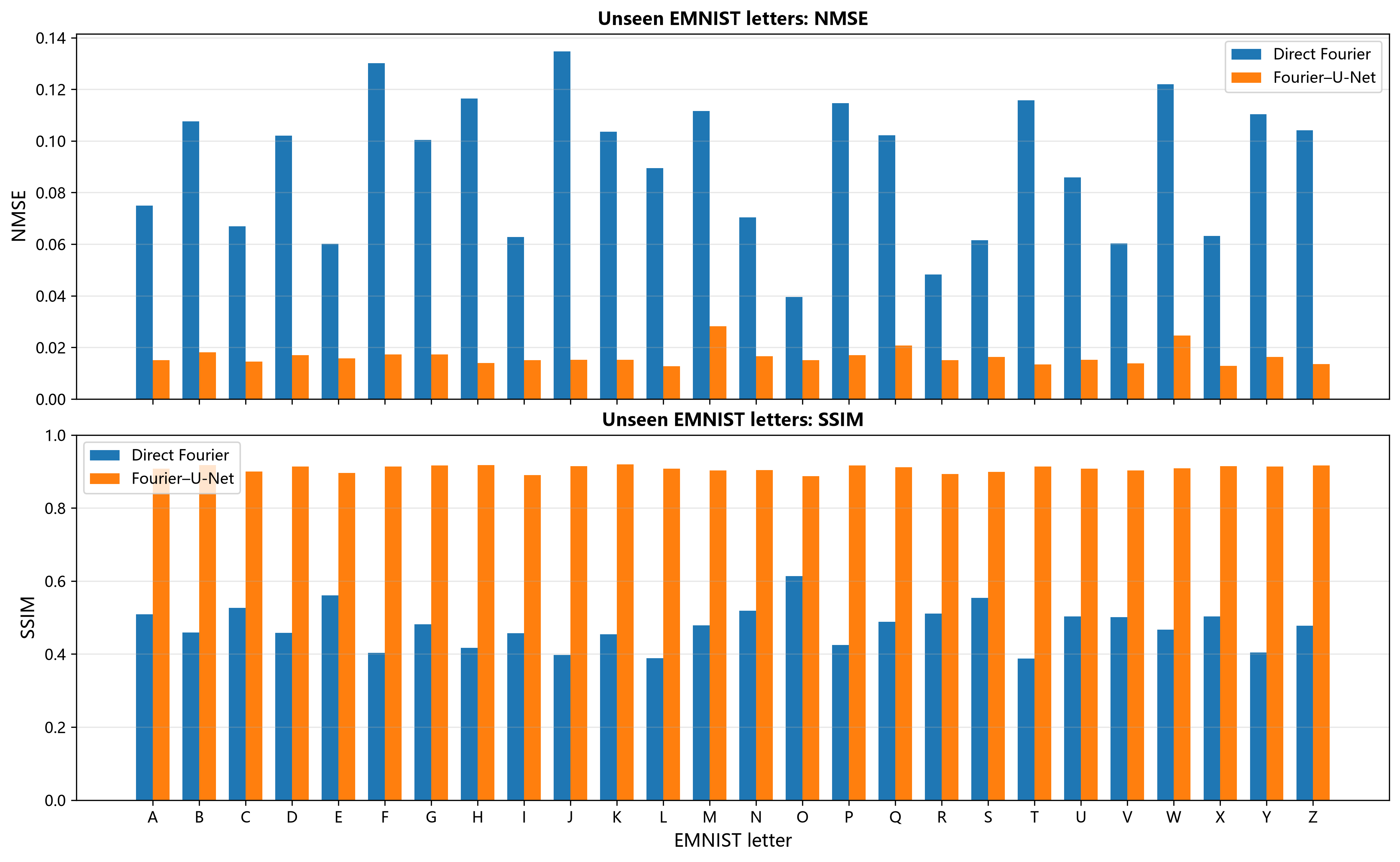}
		\caption{Comparison of NMSE and SSIM for different EMNIST letters using the frozen MNIST-trained model}
		\label{fig:exp4_letter_metric}
	\end{figure}

	Figure~\ref{fig:exp4_letter_metric} compares mean NMSE and mean SSIM for direct Fourier and Fourier--U-Net reconstructions of EMNIST letters A--Z.
	
\begin{figure}[htbp]
		\centering
		\includegraphics[
		width=0.9\textwidth,
		keepaspectratio
		]{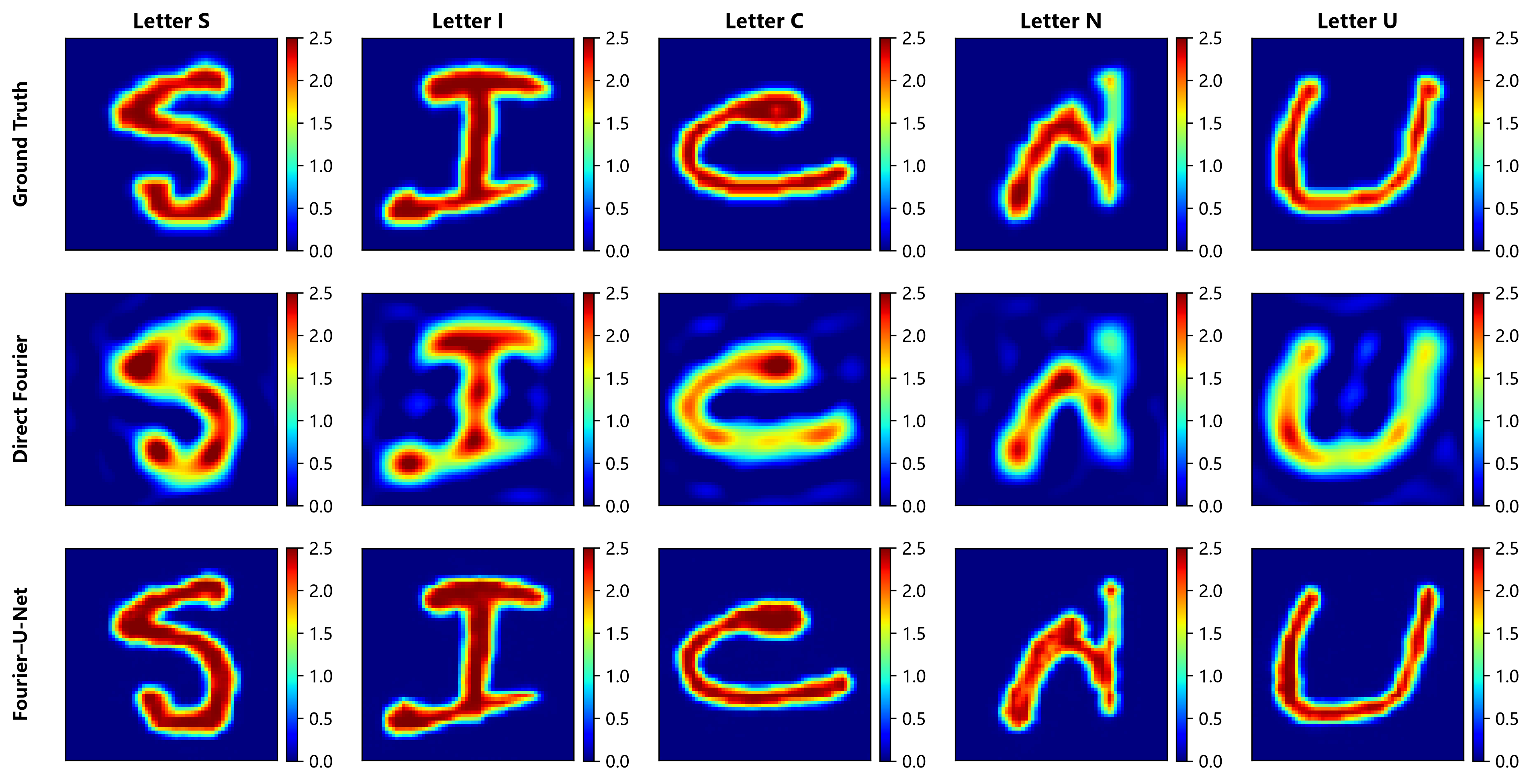}
		\caption{Reconstruction results for EMNIST handwritten letters}
		\label{fig:exp4_letter_visual}
	\end{figure}

	For every letter category in Fig.~\ref{fig:exp4_letter_metric}, Fourier--U-Net has lower mean NMSE and higher mean SSIM than direct Fourier reconstruction.
	The differences between categories may reflect variations in stroke structure and complexity.
	Figure~\ref{fig:exp4_letter_visual} shows reconstruction examples for EMNIST letters.
	
	The direct Fourier reconstructions show blurred letter strokes and local artifacts.
	The frozen U-Net better recovers the main strokes, internal gaps, and overall contours.
	
\begin{figure}[htbp]
		\centering
		\includegraphics[
		width=0.8\textwidth,
		keepaspectratio
		]{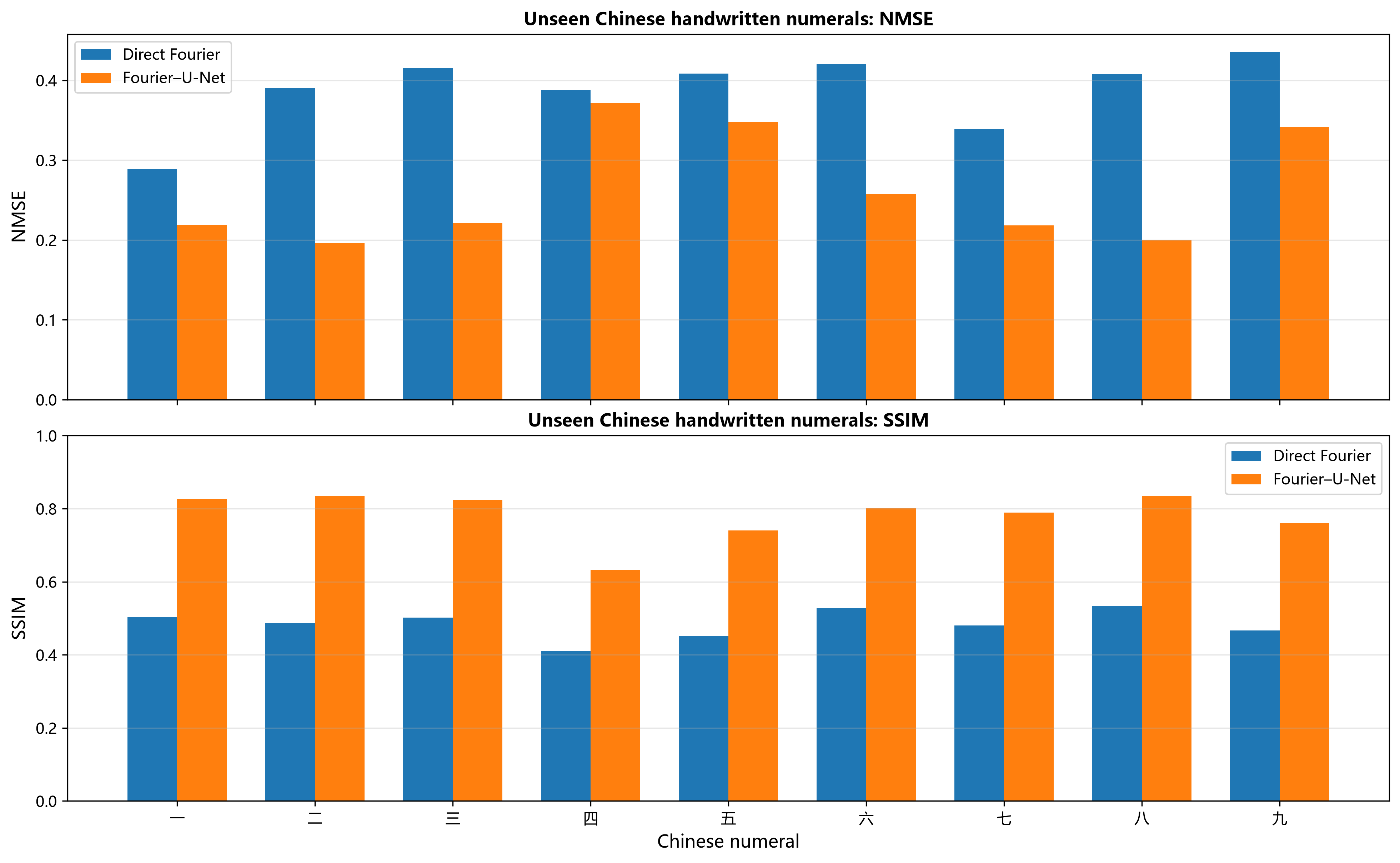}
		\caption{Comparison of NMSE and SSIM for different Chinese handwritten numerals using the MNIST-trained model}
		\label{fig:exp4_chinese_metric}
	\end{figure}

	The same frozen networks are then evaluated on Chinese handwritten numerals.
	These characters provide a further test of transfer across different stroke combinations and spatial structures.
	Figure~\ref{fig:exp4_chinese_metric} compares mean NMSE and mean SSIM for direct Fourier and Fourier--U-Net reconstructions across the nine Chinese numeral categories.
	
\begin{figure}[htbp]
		\centering
		\includegraphics[
		width=0.9\textwidth,
		keepaspectratio
		]{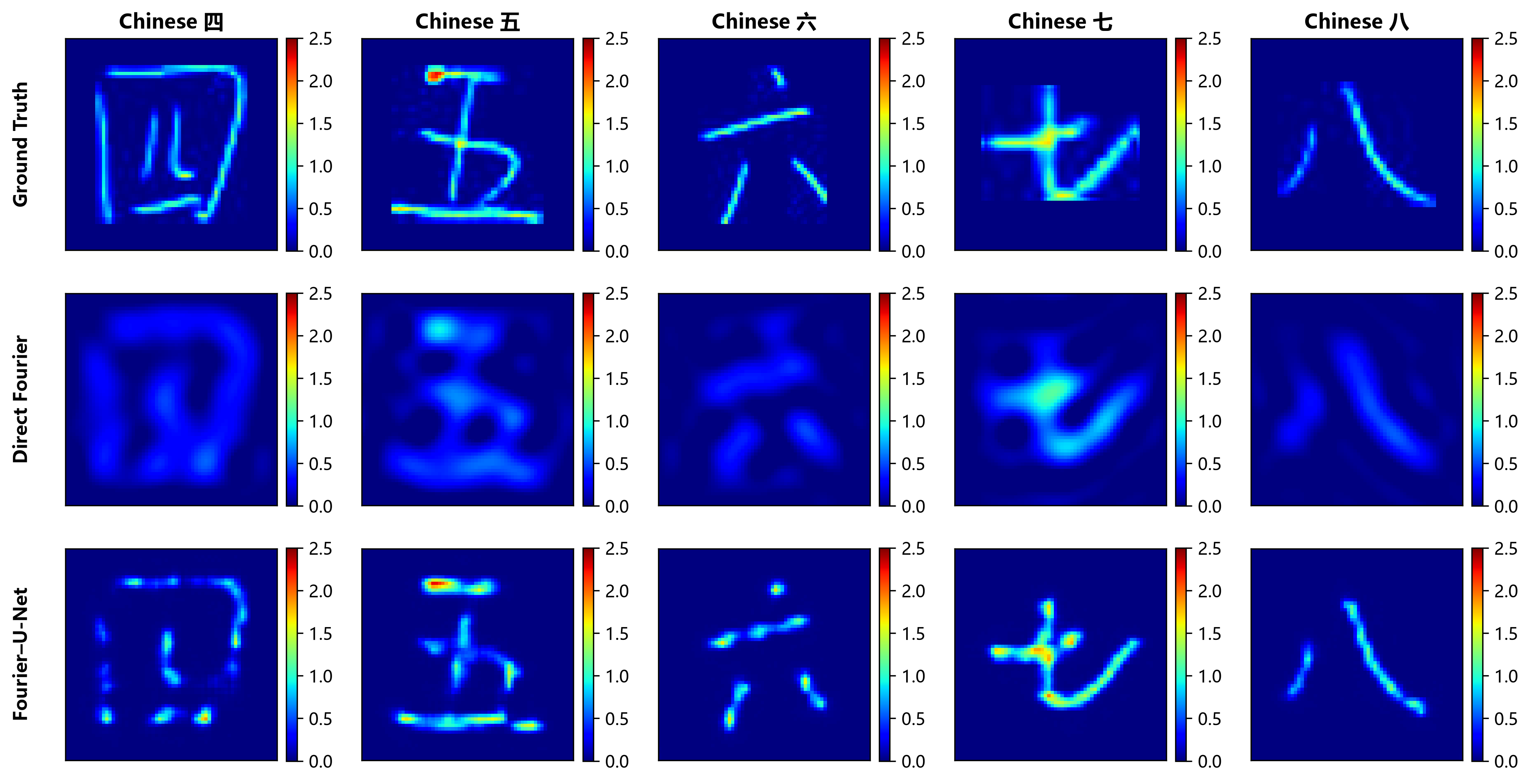}
		\caption{Reconstruction results for Chinese handwritten numeral characters}
		\label{fig:exp4_chinese_visual}
	\end{figure}

	For each Chinese numeral category in Fig.~\ref{fig:exp4_chinese_metric}, Fourier--U-Net yields lower mean NMSE and higher mean SSIM than direct Fourier reconstruction.
	The degree of improvement varies between categories.
	The Chinese numeral testing set has higher reconstruction NMSE than the EMNIST letter set.
	This difference may reflect variations in character structure relative to the MNIST training data.
	Reconstruction examples are shown in Fig.~\ref{fig:exp4_chinese_visual}.
	
	The frozen U-Net recovers the main strokes and overall structures in the displayed Chinese numeral examples.
	Some characters still show shape deviations and loss of local detail.
	Mean NMSE decreases from $0.3880$ to $0.2635$, a reduction of approximately $32.1\%$, while mean SSIM increases from $0.4847$ to $0.7828$.
	The relative reduction in NMSE and the gain in SSIM are smaller than those for EMNIST letters.
	The frozen MNIST-trained networks therefore show transfer to both tested character datasets under the same measurement conditions, improving the direct Fourier reconstructions without fine-tuning.

	\section{Conclusion}
	This paper presents a Fourier--U-Net reconstruction method for inverse scattering imaging.
	The method relates far-field scattering data to the Fourier coefficients of the target contrast under the Born approximation.
	A finite set of these coefficients is used to construct an initial contrast estimate.
	After normalization, the real and imaginary parts of this estimate form the two input channels of a U-Net that refines the reconstruction.
	
	The numerical experiments examine parameter selection, noise adaptability, complex target reconstruction, and transfer to unseen character categories.
	Experiment 1 examines how receiver distance and Fourier truncation order affect far-field extraction error, reconstruction performance, and the number of Fourier modes.
	The selected configuration $R=3$ and $N=4$ is used in the subsequent experiments.
	The comparison uses synthetic scattering data generated under the Born approximation at finite receiver distance and $5\%$ noise.
	Under these conditions, $N=4$ Fourier--U-Net has lower NMSE than direct Fourier reconstruction with either $N=4$ or $N=10$.
	The zero-frequency ablation shows that the compensation term has a substantial effect on direct Fourier reconstruction.
	After removing this term and retraining, Fourier--U-Net still has lower NMSE than direct Fourier reconstruction without compensation.
	Experiment 2 shows that high-noise pretraining followed by target-noise fine-tuning achieves reconstruction accuracy close to training from scratch under the fixed relative noise model.
	For the tested target-noise tasks and training schedules, sharing the pretrained model reduces the total training time.
	Experiment 3 shows improved direct Fourier reconstructions for MNIST single-digit and four-digit combination targets.
	Experiment 4 applies the frozen models trained on single MNIST digits to EMNIST letters and Chinese handwritten numerals without fine-tuning.
	The improvements on these unseen characters indicate transfer under the same measurement model and contrast scale.
	
	Both training and testing data are generated under the Born approximation.
	Robustness to model mismatch with full-wave nonlinear scattering data or under strong scattering conditions therefore remains untested.
	The zero-frequency compensation is constructed {\mgq by the} far-field value at a small wave number and is not extracted {\mgq from measurements}.
	The availability of this {\mgq zero-frequency} information in practical measurements requires further investigation.
	The experiments use two-dimensional simulated scattering data.
	The effects of system errors, background uncertainty, and differences between simulated and measured data on reconstruction have not yet been assessed.
	Future work will investigate reconstruction under model mismatch using full-wave nonlinear forward models and measured data, and consider extensions to three-dimensional inverse scattering imaging.
	\section*{Acknowledgments} 
	The work of Guanqiu Ma is supported by National Natural Science Foundation of China grant 12601815.

	%\section*{致谢}


\begin{thebibliography}{99}
		
		\bibitem{colton2019inverse}
		D. Colton and R. Kress,
		\textit{Inverse Acoustic and Electromagnetic Scattering Theory},
		4th ed., Applied Mathematical Sciences, vol. 93.
		Cham: Springer, 2019.
		
		\bibitem{bao2015multifrequency}
		G. Bao, P. Li, J. Lin, and F. Triki,
		``Inverse scattering problems with multi-frequencies,''
		\textit{Inverse Problems},
		vol. 31, no. 9, Art. no. 093001, 2015.
		
		\bibitem{vandenberg1997contrast}
		P. M. van den Berg and R. E. Kleinman,
		``A contrast source inversion method,''
		\textit{Inverse Problems},
		vol. 13, no. 6, pp. 1607--1620, 1997.
		
		\bibitem{potthast2006survey}
		R. Potthast,
		``A survey on sampling and probe methods for inverse problems,''
		\textit{Inverse Problems},
		vol. 22, no. 2, pp. R1--R47, 2006.
		
		\bibitem{ito2012sampling}
		K. Ito, B. Jin, and J. Zou,
		``A direct sampling method to an inverse medium scattering problem,''
		\textit{Inverse Problems},
		vol. 28, no. 2, Art. no. 025003, 2012.
		
		\bibitem{ito2013electromagnetic}
		K. Ito, B. Jin, and J. Zou,
		``A direct sampling method for inverse electromagnetic medium scattering,''
		\textit{Inverse Problems},
		vol. 29, no. 9, Art. no. 095018, 2013.
		
		\bibitem{li2013farfield}
		J. Li and J. Zou,
		``A direct sampling method for inverse scattering using far-field data,''
		\textit{Inverse Problems and Imaging},
		vol. 7, no. 3, pp. 757--775, 2013.
		
		\bibitem{zhang2015fourier}
		D. Zhang and Y. Guo,
		``Fourier method for solving the multi-frequency inverse source problem for the Helmholtz equation,''
		\textit{Inverse Problems},
		vol. 31, no. 3, Art. no. 035007, 2015.
		
		\bibitem{wang2017fourier}
		X. Wang, Y. Guo, D. Zhang, and H. Liu,
		``Fourier method for recovering acoustic sources from multi-frequency far-field data,''
		\textit{Inverse Problems},
		vol. 33, no. 3, Art. no. 035001, 2017.
		
		\bibitem{griesmaier2012windowed}
		R. Griesmaier, M. Hanke, and T. Raasch,
		``Inverse source problems for the Helmholtz equation and the windowed Fourier transform,''
		\textit{SIAM Journal on Scientific Computing},
		vol. 34, no. 3, pp. A1544--A1562, 2012.
		
		\bibitem{cheng2016stability}
		J. Cheng, V. Isakov, and S. Lu,
		``Increasing stability in the inverse source problem with many frequencies,''
		\textit{Journal of Differential Equations},
		vol. 260, no. 5, pp. 4786--4804, 2016.
		
		\bibitem{devaney1982filtered}
		A. J. Devaney,
		``A filtered backpropagation algorithm for diffraction tomography,''
		\textit{Ultrasonic Imaging},
		vol. 4, no. 4, pp. 336--350, 1982.
		
		\bibitem{chen2020review}
		X. Chen, Z. Wei, M. Li, and P. Rocca,
		``A review of deep learning approaches for inverse scattering problems (Invited Review),''
		\textit{Progress In Electromagnetics Research},
		vol. 167, pp. 67--81, 2020.
		
		\bibitem{jin2017cnn}
		K. H. Jin, M. T. McCann, E. Froustey, and M. Unser,
		``Deep convolutional neural network for inverse problems in imaging,''
		\textit{IEEE Transactions on Image Processing},
		vol. 26, no. 9, pp. 4509--4522, 2017.
		
		\bibitem{mccann2017review}
		M. T. McCann, K. H. Jin, and M. Unser,
		``Convolutional neural networks for inverse problems in imaging: A review,''
		\textit{IEEE Signal Processing Magazine},
		vol. 34, no. 6, pp. 85--95, 2017.
		
		\bibitem{arridge2019data}
		S. Arridge, P. Maass, O. \"{O}ktem, and C.-B. Sch\"{o}nlieb,
		``Solving inverse problems using data-driven models,''
		\textit{Acta Numerica},
		vol. 28, pp. 1--174, 2019.
		
		\bibitem{wei2019deep}
		Z. Wei and X. Chen,
		``Deep-learning schemes for full-wave nonlinear inverse scattering problems,''
		\textit{IEEE Transactions on Geoscience and Remote Sensing},
		vol. 57, no. 4, pp. 1849--1860, 2019.
		
		\bibitem{yao2019twostep}
		H. M. Yao, W. E. I. Sha, and L. Jiang,
		``Two-step enhanced deep learning approach for electromagnetic inverse scattering problems,''
		\textit{IEEE Antennas and Wireless Propagation Letters},
		vol. 18, no. 11, pp. 2254--2258, 2019.
		
		\bibitem{zhou2021improved}
		Y. Zhou, Y. Zhong, Z. Wei, T. Yin, and X. Chen,
		``An improved deep learning scheme for solving 2-D and 3-D inverse scattering problems,''
		\textit{IEEE Transactions on Antennas and Propagation},
		vol. 69, no. 5, pp. 2853--2863, 2021.
		
		\bibitem{gao2022ann}
		Y. Gao, H. Liu, X. Wang, and K. Zhang,
		``On an artificial neural network for inverse scattering problems,''
		\textit{Journal of Computational Physics},
		vol. 448, Art. no. 110771, 2022.
		
		\bibitem{ronneberger2015unet}
		O. Ronneberger, P. Fischer, and T. Brox,
		``U-Net: Convolutional networks for biomedical image segmentation,''
		in \textit{Medical Image Computing and Computer-Assisted Intervention--MICCAI 2015},
		Lecture Notes in Computer Science, vol. 9351,
		pp. 234--241, Springer, Cham, 2015.
		
		\bibitem{benny2022microwave}
		R. Benny, T. A. Anjit, and P. Mythili,
		``Deep learning based non-iterative solution to the inverse problem in microwave imaging,''
		\textit{Progress In Electromagnetics Research M},
		vol. 109, pp. 231--240, 2022.
		
		\bibitem{zhang2020hybrid}
		L. Zhang, K. Xu, R. Song, X. Ye, G. Wang, and X. Chen,
		``Learning-based quantitative microwave imaging with a hybrid input scheme,''
		\textit{IEEE Sensors Journal},
		vol. 20, no. 24, pp. 15007--15013, 2020.
		
		\bibitem{ning2024dsm}
		J. Ning, F. Han, and J. Zou,
		``A direct sampling-based deep learning approach for inverse medium scattering problems,''
		\textit{Inverse Problems},
		vol. 40, no. 1, Art. no. 015005, 2024.
		
		\bibitem{chen2026fourier}
		H. Chen, Y. Chang, Y. Guo, and Y. Wang,
		``A deep learning-enhanced Fourier method for the multi-frequency inverse source problem with sparse far-field data,''
		\textit{Communications in Computational Physics},
		vol. 40, no. 5, pp. 1592--1610, 2026.
		
		\bibitem{wang2004ssim}
		Z. Wang, A. C. Bovik, H. R. Sheikh, and E. P. Simoncelli,
		``Image quality assessment: From error visibility to structural similarity,''
		\textit{IEEE Transactions on Image Processing},
		vol. 13, no. 4, pp. 600--612, 2004.
		
		\bibitem{lecun1998gradient}
		Y. LeCun, L. Bottou, Y. Bengio, and P. Haffner,
		``Gradient-based learning applied to document recognition,''
		\textit{Proceedings of the IEEE},
		vol. 86, no. 11, pp. 2278--2324, 1998.
		
		\bibitem{cohen2017emnist}
		G. Cohen, S. Afshar, J. Tapson, and A. van Schaik,
		``EMNIST: Extending MNIST to handwritten letters,''
		in \textit{2017 International Joint Conference on Neural Networks (IJCNN)},
		pp. 2921--2926, 2017.
		
		\bibitem{nazarpour2017chinese}
		K. Nazarpour and M. Chen,
		\textit{Handwritten Chinese Numbers},
		Newcastle University, Dataset, 2017.
		
	\end{thebibliography}
\end{document}